\documentclass[preprint]{elsarticle}
\usepackage{amssymb}
\usepackage{amsmath}
\usepackage{amsthm}
\newcommand{\R}{{\Bbb R}}
\newcommand{\C}{{\Bbb C}}
\newcommand{\N}{{\Bbb N}}
\newcommand{\D}{{\cal D}}
\newcommand{\A}{{\cal A}}
\renewcommand{\S}{{\cal S}}

\newcommand{\cL}{\mathcal{L}}
\newcommand{\M}{{\rm M}}

\newcommand{\Var}{{\rm Var\,}}
\newcommand{\const}{{\rm const}}
\newcommand{\supp}{\mathop{\rm supp}}

\newcommand{\pr}{{\rm pr}}
\renewcommand{\limsup}{\mathop{\overline{\lim}}}
\renewcommand{\liminf}{\mathop{\underline{\lim}}}

\renewcommand{\Im}{\mathop{\rm Im}}
\newcommand{\esslim}{\mathop{\rm ess\,lim}}
\newcommand{\esssup}{\mathop{\rm ess\,sup}}
\newcommand{\esslimsup}{\mathop{\rm ess\,limsup}}
\newcommand{\Blim}{\mathop{\rm B\,lim}}

\newcommand{\<}{\langle}
\renewcommand{\>}{\rangle}
\newtheorem{thm}{Theorem}[section]
\newtheorem{lem}{Lemma}[section]
\newtheorem{prop}{Proposition}[section]
\newtheorem{cor}{Corollary}[section]
\newdefinition{defn}{Definition}
\newdefinition{rem}{Remark}[section]
\newproof{pf}{Proof}
\makeatletter \@addtoreset{equation}{section}

\makeatother

\journal{ArXiv}

\begin{document}
\begin{frontmatter}

\title{On variants of $H$-measures and compensated compactness
\tnoteref{grants}}
\tnotetext[grants]{This research was carried out under financial support of the Russian Foundation for
Basic Research (grant No.~12-01-00230-a).}
\author{E.Yu.~Panov}
\ead{Eugeny.Panov@novsu.ru}


\begin{abstract}
We introduce new variant of $H$-measures defined on spectra of general algebra of test symbols and derive the localization properties of such $H$-measures.
Applications for the compensated compactness theory are given. In particular, we present new compensated compactness results for quadratic functionals in the case of
general pseudo-differential constraints. The case of inhomogeneous second order differential constraints is also studied.
\end{abstract}

\end{frontmatter}

\section{Introduction}\label{Se1}

Let $$F(u)(\xi)=\int_{\R^n} e^{-2\pi i\xi\cdot x} u(x)dx, \quad \xi\in\R^n,$$
be the Fourier transformation extended as a unitary operator on the space $u(x)\in L^2(\R^n)$,
let $S=S^{n-1}= \{ \ \xi\in\R \ \mid \ |\xi|=1 \ \}$
be the unit sphere in $\R^n$.
Denote by $u\to\overline{u}$, $u\in\C$ the complex conjugation.

 The  concept  of an $H$-measure corresponding to some
 sequence  of  vector-valued
  functions bounded in $L^2(\Omega)$, where $\Omega\subset\R^n$ is an open domain, was introduced by Tartar \cite{Ta2} and Ger\'ard \cite{Ger}
on the  basis of  the following result. For $r\in\N$ let
$U_r(x)=\left(U_r^1(x),\ldots,U_r^N(x)\right)\in L^2(\Omega, \R^N)$
 be a sequence weakly convergent to
the zero vector.

\begin{prop}[see Theorem~1.1 in \cite{Ta2}]\label{Pro1} There
exists a family of complex Borel measures
$\mu=\left\{\mu^{\alpha\beta}\right\}_{\alpha,\beta=1}^N$
 in $\Omega\times  S$ and a subsequence of $U_r(x)$ (still denoted $U_r$) such that
      \begin{equation}\label{Hm}
\langle\mu^{\alpha\beta},\Phi_1(x)\overline{\Phi_2(x)}\psi(\xi)\rangle=
\lim\limits_{r\to\infty}\int_{\R^n}
F(U_r^\alpha\Phi_1)(\xi)\overline{F(U_r^\beta\Phi_2)(\xi)}
\psi\left(\frac{\xi}{|\xi|}\right)d\xi
   \end{equation}
for all $\Phi_1(x),\Phi_2(x)\in C_0(\Omega)$ and $\psi(\xi)\in
C(S)$.
\end{prop}
Here and in the sequel we use notations $C_0(\Omega)$ for the space of continuous functions on $\Omega$ with compact supports.

The family $\mu=\left\{\mu^{\alpha\beta}\right\}_{\alpha,\beta=1}^N$ is called the
$H$-measure corresponding to $U_r(x).$

In \cite{AntLaz} the new concept of parabolic $H$-measures was suggested. This concept was extended in \cite{PaJMS}, where the notion of ultra-parabolic $H$-measures was introduced. Suppose that $X\subset\R^n$ is a linear subspace, $X^\bot$
is its orthogonal complement, $P_1,P_2$ are orthogonal projections on $X$, $X^\bot$, respectively. We
denote for $\xi\in\R^n$ $\tilde\xi=P_1\xi$, $\bar\xi=P_2\xi$, so that $\tilde\xi\in X$, $\bar\xi\in
X^\bot$, $\xi=\tilde\xi+\bar\xi$. Let $S_X=\{ \ \xi\in\R^n \ | \ |\tilde\xi|^2+|\bar\xi|^4=1 \ \}$.
Then $S_X$ is a compact smooth manifold of codimension $1$; in the case when $X=\{0\}$ or $X=\R^n$, it
coincides with the unit sphere $S=\{ \xi\in\R^n \ | \ |\xi|=1 \ \}$. Let us define a projection
$\pi_X:\R^n\setminus\{0\}\to S_X$ by
$$
\pi_X(\xi)=\frac{\tilde\xi}{(|\tilde\xi|^2+|\bar\xi|^4)^{1/2}}+\frac{\bar\xi}{(|\tilde\xi|^2+|\bar\xi|^4)^{1/4}}.
$$
Remark that in the case when $X=\{0\}$ or $X=\R^n$, \ $\pi_X(\xi)=\xi/|\xi|$ is the orthogonal projection on the sphere.
With the notations from Proposition~\ref{Pro1}, the following extension holds:

\medskip
\begin{prop}[see \cite{PaJMS,AIHP}]\label{Pro2} There exists a family of complex Borel measures
$\mu=\left\{\mu^{\alpha\beta}\right\}_{\alpha,\beta=1}^N$
 in $\Omega\times  S_X$ and a subsequence $U_r(x)=U_k(x),$ $k=k_r,$ such that
\begin{equation}\label{Hm'}
\<\mu^{\alpha\beta},\Phi_1(x)\overline{\Phi_2(x)}\psi(\xi)\>=
\lim\limits_{r\to\infty}\int_{\R^n}
F(U_r^\alpha\Phi_1)(\xi)\overline{F(U_r^\beta\Phi_2)(\xi)}
\psi\left(\pi_X(\xi)\right)d\xi
\end{equation}
for all $\Phi_1(x),\Phi_2(x)\in C_0(\Omega)$ and $\psi(\xi)\in
C(S_X)$.
\end{prop}

The family $\mu=\left\{\mu^{\alpha\beta}\right\}_{\alpha,\beta=1}^N$ we shall call an ultra-parabolic $H$-measure corresponding to $U_r(x).$

In paper \cite{AIHP} the localization properties of ultra-parabolic $H$-measures were applied to extend
the compensated compactness theory \cite{Mu,Ta1} for weakly convergent sequences $u_r\in L^p_{loc}(\Omega,\R^N)$ to the case when the differential constraints may contain second-order terms while all the coefficients are variable. We describe the results of \cite{AIHP} in the particular case $p=2$. Thus, assume that
a sequence $u_r\in L^2_{loc}(\Omega,\R^N)$ converges weakly to a vector-function $u(x)$ as $r\to\infty$
and satisfies the condition that
the sequences
\begin{equation}
\label{HP1} \sum_{\alpha=1}^N\sum_{k=1}^n \partial_{x_k}(a_{s\alpha k}u_{\alpha r})
+\sum_{\alpha=1}^N\sum_{k,l=\nu+1}^n \partial_{x_kx_l}(b_{s\alpha kl}u_{\alpha r}), \quad s=1,\ldots,m
\end{equation}
are pre-compact in the anisotropic Sobolev space $W^{-1,-2}_{2,loc}(\Omega)$ (~the parameter $-1$ corresponds to the first $\nu$ variables $x_1,\ldots, x_\nu$ while the parameter $-2$ corresponds to the remaining variables $x_{\nu+1},\ldots, x_n$~). Here $\nu$ is an integer number between $0$ and $n$, and the coefficients
$a_{s\alpha k}=a_{s\alpha k}(x)$, $b_{s\alpha kl}=b_{s\alpha kl}(x)$ are assumed to be continuous on $\Omega$.

We introduce the set $\Lambda$ (here $i=\sqrt{-1}$):
\begin{eqnarray}
\label{HP2}
\Lambda=\Lambda(x)=\Bigl\{ \lambda\in\C^N \ | \ \exists\xi\in\R^n, \xi\not=0: \nonumber \\
\sum_{\alpha=1}^N\Bigl(i\sum_{k=1}^\nu a_{s\alpha k}(x)\xi_k-\sum_{k,l=\nu+1}^n b_{s\alpha
kl}(x)\xi_k\xi_l\Bigr)\lambda_\alpha=0 \ \forall s=1,\ldots,m \ \Bigr\}.
\end{eqnarray}
Consider the quadratic form $q(x,u)=Q(x)u\cdot u$, where $Q(x)$ is a symmetric matrix with
coefficients $q_{\alpha\beta}(x)\in C(\Omega)$, $\alpha,\beta=1,\ldots,N$ and $u\cdot v$ denotes the scalar
multiplication on $\R^N$. The form $q(x,u)$ can be extended as Hermitian form on $\C^N$ by the
standard relation
$$
q(x,u)=\sum_{\alpha,\beta=1}^N q_{\alpha\beta}(x)u_\alpha \overline{u_\beta}.
$$
Now, let the sequence $q(x,u_r)\rightharpoonup v$ as $r\to\infty$ weakly in
$\D'(\Omega)$. Since this
sequence is bounded in
$L^1_{loc}(\Omega)$ then,
passing to a subsequence if necessary, we may claim that $v$ is a locally finite measure on $\Omega$ (~i.e., $v\in M_{loc}(\Omega)$~), and $q(x,u_r)\rightharpoonup v$ weakly in $M_{loc}(\Omega)$.
The following result was established in \cite{AIHP}.

\begin{thm}\label{Th1} Assume that $q(x,\lambda)\ge 0$ for all $\lambda\in\Lambda(x)$,
$x\in\Omega$. Then ${q(x,u(x))\le v}$ ( in the sense of measures ).
\end{thm}

In the case $\nu=n$ when the second order terms in (\ref{HP1}) are absent and all the coefficients are constant the statement of Theorem~\ref{Th1} is the classical Tartar-Murat compensated compactness.

In this paper we generalize the result of Theorem~\ref{Th1} to the case when the degeneration subspaces $X_s$ in constraints (\ref{HP1}) may depend on $s$ and give some applications.

For that, we introduce the general variant of $H$-measures by extension of a class of admissible test functions $\psi(\xi)$. We will describe this class in the next section.

\section{Algebra of admissible symbols}\label{Se2}

Let us denote by $B_\Phi$ and
$A_\psi$ the bounded pseudodifferential operators on $L^2(\R^n)$ with symbols $\Phi(x),\psi(\xi)\in L^\infty(\R^n)$, respectively, that is,
$$
B_\Phi u(x)=\Phi(x)u(x), \quad F (A_\psi u)(\xi)=\psi(\xi) F(u)(\xi).
$$
We introduce the subalgebra $A$ of the algebra $L^\infty(\R^n)$, consisting of bounded measurable functions $\psi(\xi)$ on $\R^n$ such that the commutators $[A_\psi,B_\Phi]$ are
compact operators in $L^2(\R^n)$ for all $\Phi(x)\in C_0(\R^n)$. Let $A_0=L^\infty_0$ be a subspace of $L^\infty(\R^n)$ consisting of functions $\psi(\xi)$ vanishing at infinity: $\esslim\limits_{|\xi|\to\infty}\psi(\xi)=0$.

\begin{lem}\label{Lem1} For every $\Phi(x)\in C_0(\R^n)$, $\psi(\xi)\in A_0$ the operators
$A_\psi B_\Phi$, $B_\Phi A_\psi$ are compact in $L^2(\R^n)$.
\end{lem}

\begin{pf}
First, assume that $\psi(\xi)\in L^\infty(\R^n)$ is a function with compact support $K=\supp\psi\subset\R^n$.
Let $u_k$, $k\in\N$, be a sequence in $L^2(\R^n)$, weakly convergent to zero: $u_k\mathop{\rightharpoonup}\limits_{k\to\infty} 0$. We have to prove that
$A_\psi B_\Phi u_k\mathop{\to}\limits_{k\to\infty} 0$ in $L^2(\R^n)$ (strongly).
Since $B_\Phi u_k=\Phi(x)u_k(x)\mathop{\rightharpoonup}\limits_{k\to\infty} 0$ weakly in $L^1(\R^n)$, then
$$
F(B_\Phi u_k)(\xi)=\int_{\R^n}e^{-2\pi i\xi\cdot x}u_k(x)\Phi(x)dx\mathop{\to}_{k\to\infty} 0
$$
for all $\xi\in\R^n$, and
$$
|F(B_\Phi u_k)(\xi)|\le \|\Phi u_k\|_1\le C=\|\Phi\|_2\sup\limits_{k\in\N}\|u_k\|_2<\infty.
$$
Then, by the Lebesgue dominated convergence theorem, we claim that
$$
\|A_\psi B_\Phi u_k\|^2_2=\int_K |F(B_\Phi u_k)(\xi)\psi(\xi)|^2d\xi\to 0
$$
as $k\to \infty$, that is, $A_\psi B_\Phi u_k\to 0$ in $L^2(\R^n)$. We see that the operator $A_\psi B_\Phi$ transforms weakly convergent sequences in $L^2$ to strongly convergent ones. Hence, this operator is compact.

In the general case $\psi(\xi)\in A_0$ we introduce the sequence $\psi_m(\xi)=\psi(\xi)\theta(m-|\xi|)$, $m\in\N$,
where $\theta(r)=\left\{\begin{array}{lr} 0, & r\le 0, \\ 1, & r>0 \end{array}\right.$ is the Heaviside function.
Then
$$\|\psi_m-\psi\|_\infty=\esssup\limits_{|\xi|\ge m}|\psi(\xi)|\to 0$$
as $m\to\infty$, and therefore the operator norms
$$
\|A_{\psi_m}-A_\psi\|=\|\psi_m-\psi\|_\infty\mathop{\to}_{m\to\infty} 0.
$$
This implies that $A_{\psi_m}B_\Phi\to A_\psi B_\Phi$ as $m\to\infty$ in the algebra of bounded linear operators on
$L^2(\R^n)$. The functions $\psi_m(\xi)$ have compact supports and it has been already proven that the operators $A_{\psi_m}B_\Phi$ are compact.
We conclude that $A_\psi B_\Phi$ is a compact operator, as the limit of the sequence of compact operators $A_{\psi_m}B_\Phi$.

In order to prove compactness of $B_\Phi A_\psi$, observe that this operator is conjugate to
$A_{\bar\psi} B_{\bar\Phi}= (A_\psi)^* (B_\Phi)^*$. As we have already established, the operator $A_{\bar\psi} B_{\bar\Phi}$ is compact. Therefore, the operator $B_\Phi A_\psi=(A_{\bar\psi} B_{\bar\Phi})^*$ is compact as well. The proof is complete.
\end{pf}

In view of Lemma~\ref{Lem1} we find that for $\psi(\xi)\in A_0$ the commutator $[A_\psi,B_\Phi]=A_\psi B_\Phi-B_\Phi A_\psi$
is a compact operator in $L^2(\R^n)$ for all $\Phi(x)\in C_0(\R^n)$. In particular $A_0\subset A$. It is clear that  $A_0$ is a closed ideal in $A$. We denote by $\A=A/A_0$ the correspondent quotient
algebra. Clearly, $\A$ is a commutative  Banach $C^*$-algebra (subject to the involution defined by
complex conjugation) equipped with the factor-norm (we identify the class $[\psi]\in\A$ with the corresponding representative function $\psi(\xi)$)
$$\|\psi\|=\esslimsup_{\xi\to\infty}|\psi(\xi)|=\lim\limits_{R\to\infty}\esssup\limits_{|\xi|>R}|\psi(\xi)|.$$
Therefore, the Gelfand transform $\psi(\xi)\to\hat\psi(\eta)$ is an
isomorphism of $\A$ into the algebra $C(\S)$ of continuous functions on the spectrum $\S$ of $\A$.

We introduce the order in $\A$ generated by the cone of nonnegative functions, that is,
a class $a\ge 0$ if and only if there exists a real nonnegative function $\psi\in a$, i.e., $a=[\psi]$. As is easy to verify, for $a,b\in\A$, $a,b\ge 0$, and $\alpha,\beta\in [0,+\infty)$ \
$\alpha a+\beta b\ge 0$ $ab\ge 0$. As usual, we say that $a_1\ge a_2$ if $a_1-a_2\ge 0$. It turns out that the Gelfand transform is monotone, that is, the following statement is fulfilled.

\begin{lem}\label{Lem2} The class $a=[\psi]\ge 0$ if and only if $\hat\psi(\eta)\ge 0$ for all $\eta\in\S$.
\end{lem}
\begin{pf}
If  $\hat\psi(\eta)\ge 0$ for all $\eta\in\S$ then the function $\alpha(\eta)=(\hat\psi(\eta))^{1/2}$ is well-defined and continuous on $\S$. Therefore, there exists a unique class $b=[\beta(\xi)]\in\A$ such that $\alpha(\eta)=\hat\beta(\eta)$. Since the Gelfand transform satisfies the property $\widehat{\bar\psi}(\eta)=\overline{\hat\psi}(\eta)$, we see that $\widehat{b\bar b}(\eta)=
(\alpha(\eta))^2=\hat\psi(\eta)$ and the equality $a=[\psi]=b\bar b=[|\beta|^2]$ follows. This equality implies that $a\ge 0$.

Conversely, let $a=[\psi]\ge 0$. Since $a=\bar a$, the function $\hat\psi(\eta)$ is real. We define the real nonnegative functions $\hat\psi^\pm(\eta)=\max(0,\pm\hat\psi(\eta))\in C(\S)$. Then, there exist
classes $a^\pm=[\psi^\pm]$ such that $\widehat{\psi^\pm}(\eta)=\hat\psi^\pm(\eta)$. As we have already established, $a^\pm\ge 0$. Since $\hat\psi(\eta)=\hat\psi^+(\eta)-\hat\psi^-(\eta)$, and $\hat\psi^+(\eta)\cdot\hat\psi^-(\eta)=0$, the same is true for the $a^\pm$:
$a=a^+-a^-$, $a^+a^-=0$. Therefore, $-aa^-=(a^-)^2\ge 0$. On the other hand, $aa^-\ge 0$, as a product of nonnegative elements. We conclude that $(a^-)^2=-aa^-=0$ and, therefore, $a^-=0$. But this means that $\hat\psi^-(\eta)=0$ and implies nonnegativity of $\hat\psi(\eta)$: $\hat\psi(\eta)=\hat\psi^+(\eta)\ge 0$.
This completes the proof.
\end{pf}

As follows from \cite[Lemma~2]{PaJMS}, functions $\psi(\pi_X(\xi))$ belong to the algebra $A$ for each
$\psi\in C(S_X)$. Hence, the algebra of quasi-homogeneous functions $$A_X=\{ \ \psi(\pi_X(\xi)) \ | \ \psi\in C(S_X)
\ \}$$ is a closed $C^*$-subalgebra of $\A$ and its spectrum coincides with $S_X$. The embedding
$A_X\subset\A$ yields the continuous projection of the spectra $p_X:\S\to S_X$.  One of our aims is to formulate localization properties for $H$-measures corresponding to
sequences satisfying general second order differential constraints. For this, we need to find simple necessary and sufficient
conditions for a family of vectors $\{\xi_X\}_{X\subset\R^n}$ to satisfy the property
$\xi_X=p_X(\eta)$ for all $X\subset\R^n$, where $\eta\in\S$.
The following statement holds.

\begin{prop}\label{Pro3} Assume that $\eta\in\S$ and for $X\subset\R^n$ let $p_X(\eta)=(\tilde\xi_X,\bar\xi_X)\in X\oplus X^\bot$. Then there exist a unique orthonormal system $\{\zeta_1,\ldots,\zeta_m\}$ in $\R^n$ and an integer $d\in\{m-1,m\}$ such that

(i) $\tilde\xi_X\not=0\Leftrightarrow X\supset \tilde X\doteq\cL(\zeta_1,\ldots,\zeta_d)$ (this is a linear span of vectors $\zeta_1,\ldots,\zeta_d$). Besides, if $\tilde\xi_X\not=0$, then $\tilde\xi_X\upuparrows\zeta_1$;

(ii) $\bar\xi_X\not=0\Leftrightarrow X\not\supset \bar X\doteq\cL(\zeta_1,\ldots,\zeta_m)$. Besides, if $\bar\xi_X\not=0$, then $\bar\xi_X\upuparrows \pr_{X^\bot}\zeta_{k(X)}$, where
$k(X)=\min\{ \ k=1,\ldots,m \ | \ \zeta_k\notin X \ \}$.
\end{prop}
\begin{pf}
We divide the proof into 6 steps.

{\bf 1st Step.}

We introduce the set $\tilde\cL$ of all subspaces $X\subset\R^n$ such that $\tilde\xi_X\not=0$. Let us show that $\tilde\cL$ contains the smallest space. For that, we first prove that the intersection $X_1\cap X_2$ of two spaces
$X_1,X_2\in\tilde\cL$ lays in $\tilde\cL$ as well. We denote $X_0=X_1\cap X_2$, $X_{10}=X_1\ominus X_0=\{x\in X_1: x\bot X_0\}$, $X_{20}=X_2\ominus X_0$. Then we have the following representations
\begin{equation}\label{10}
\R^n=X_0\oplus X_{10}\oplus X_1^\bot=X_0\oplus X_{20}\oplus X_2^\bot.
\end{equation}
Let
\begin{equation}\label{11}\xi=\xi_0+\xi_1+\xi_3=\xi_0+\xi_2+\xi_4
\end{equation}
be orthogonal decompositions of a vector $\xi\in\R^n$ corresponding to (\ref{10}). Here $\xi_0\in X_0$, $\xi_1\in X_{10}$, $\xi_2\in X_{20}$, $\xi_3\in X_1^\bot$, and $\xi_4\in X_2^\bot$.
We introduce the functions
$$
f_1(\xi)=\frac{|\xi_0|^2+|\xi_1|^2}{|\xi_0|^2+|\xi_1|^2+|\xi_3|^4}, \quad
f_1(\xi)=\frac{|\xi_0|^2+|\xi_2|^2}{|\xi_0|^2+|\xi_2|^2+|\xi_4|^4}
$$
defined on $\R^n\setminus\{0\}$. Obviously, $f_1\in A_{X_1}\subset\A$, $f_2\in A_{X_2}\subset\A$, and
$\widehat{f_1}(\eta)=|\tilde\xi_{X_1}|^2\not=0$, $\widehat{f_2}(\eta)=|\tilde\xi_{X_2}|^2\not=0$.
We define the subspace
$Y\subset X_1^\bot\oplus X_2^\bot$ consisting of pairs $(\xi_3,\xi_4)$ such that $\xi_1+\xi_3=\xi_2+\xi_4$ for some
vectors $\xi_1\in X_{10}$, $\xi_2\in X_{20}$. Observe that the vectors $\xi_1$, $\xi_2$ are uniquely defined by the above equality. Indeed, if $\xi'_1+\xi_3=\xi'_2+\xi_4$ for some other vectors $\xi'_1\in X_{10}$, $\xi'_2\in X_{20}$
then $\xi_1-\xi'_1=\xi_2-\xi'_2\in X_{10}\cap X_{20}=\{0\}$ and we conclude that
$\xi_1=\xi'_1$, $\xi_2=\xi'_2$. Thus, we can define  the linear maps $A_1:Y\to X_{10}$, $A_2:Y\to X_{20}$ such that
$A_1(\xi_3,\xi_4)=\xi_1$, $A_2(\xi_3,\xi_4)=\xi_2$. Since these maps are continuous, we can find a positive constant
$C$ such that
\begin{equation}\label{est}
|A_i(\xi_3,\xi_4)|^2\le C(|\xi_3|^2+|\xi_4|^2) \ \mbox{ for all } (\xi_3,\xi_4)\in Y.
\end{equation}
 Then
\begin{eqnarray}\label{12}
f_1(\xi)\le\frac{|\xi_0|^2+C(|\xi_3|^2+|\xi_4|^2)}{|\xi_0|^2+|\xi_1|^2+|\xi_3|^4}\le
\frac{|\xi_0|^2+C|\xi_4|^2}{|\xi_0|^2+|\xi_1|^2+|\xi_3|^4}+\alpha_1(\xi), \nonumber\\
f_2(\xi)\le\frac{|\xi_0|^2+C(|\xi_3|^2+|\xi_4|^2)}{|\xi_0|^2+|\xi_2|^2+|\xi_4|^4}\le
\frac{|\xi_0|^2+C|\xi_3|^2}{|\xi_0|^2+|\xi_2|^2+|\xi_4|^4}+\alpha_2(\xi),
\end{eqnarray}
where
\begin{eqnarray*}
\alpha_1(\xi)=\frac{C|\xi_3|^2}{|\xi_0|^2+|\xi_1|^2+|\xi_3|^4}\mathop{\to}_{\xi\to\infty} 0, \\
\alpha_2(\xi)=\frac{C|\xi_4|^2}{|\xi_0|^2+|\xi_2|^2+|\xi_4|^4}\mathop{\to}_{\xi\to\infty} 0,
\end{eqnarray*}
that is, $\alpha_k(\xi)\in A_0$, $j=1,2$. In view of (\ref{12})
\begin{eqnarray*}
0\le f_1(\xi)\le\frac{|\xi_0|^2+C|\xi_4|^2}{|\xi_0|^2+|\xi_1|^2+|\xi_3|^4}, \\
0\le f_2(\xi)\le \frac{|\xi_0|^2+C|\xi_3|^2}{|\xi_0|^2+|\xi_2|^2+|\xi_4|^4}
\end{eqnarray*}
in $\A$, which implies that in this algebra
\begin{equation}\label{13}
0\le f_1(\xi)f_2(\xi)\le\frac{(|\xi_0|^2+C|\xi_4|^2)(|\xi_0|^2+C|\xi_3|^2)}
{(|\xi_0|^2+|\xi_1|^2+|\xi_3|^4)(|\xi_0|^2+|\xi_2|^2+|\xi_4|^4)}.
\end{equation}
Observe that
\begin{eqnarray*}
(|\xi_0|^2+C|\xi_4|^2)(|\xi_0|^2+C|\xi_3|^2)\le |\xi_0|^4+C|\xi_4|^2(|\xi_0|^2+C|\xi_3|^2)+C|\xi_0|^2|\xi_3|^2,\\
(|\xi_0|^2+|\xi_1|^2+|\xi_3|^4)(|\xi_0|^2+|\xi_2|^2+|\xi_4|^4)\ge(|\xi_0|^2+|\xi_3|^4)(|\xi_0|^2+|\xi_4|^4),
\end{eqnarray*}
and it follows from (\ref{13}) that
\begin{eqnarray}\label{14}
0\le f_1(\xi)f_2(\xi)\le\nonumber\\ \frac{|\xi_0|^4}{(|\xi_0|^2+|\xi_3|^4)(|\xi_0|^2+|\xi_4|^4)}+
C\frac{|\xi_0|^2+C|\xi_3|^2}{|\xi_0|^2+|\xi_3|^4}\cdot\frac{|\xi_4|^2}{|\xi_0|^2+|\xi_4|^4}+\nonumber\\
C\frac{|\xi_0|^2}{|\xi_0|^2+|\xi_4|^4}\cdot\frac{|\xi_3|^2}{|\xi_0|^2+|\xi_3|^4}=
\frac{|\xi_0|^4}{(|\xi_0|^2+|\xi_3|^4)(|\xi_0|^2+|\xi_4|^4)}+\beta(\xi),
\end{eqnarray}
where $\beta(\xi)\in A_0$. Since
\begin{eqnarray*}
(|\xi_1|^2+|\xi_3|^2)^2\le (C(|\xi_3|^2+|\xi_4|^2)+|\xi_3|^2)^2 \le\\ (C+1)^2(|\xi_3|^2+|\xi_4|^2)^2
\le 2(C+1)^2(|\xi_3|^4+|\xi_4|^4),
\end{eqnarray*}
then
\begin{eqnarray*}
(|\xi_0|^2+|\xi_3|^4)(|\xi_0|^2+|\xi_4|^4)\ge |\xi_0|^2(|\xi_0|^2+|\xi_3|^4+|\xi_4|^4)\ge \\
\frac{1}{2(C+1)^2}|\xi_0|^2(|\xi_0|^2+(|\xi_1|^2+|\xi_3|^2)^2),
\end{eqnarray*}
and it follows from (\ref{14}) that in $\A$
\begin{equation}\label{15}
0\le f_1(\xi)f_2(\xi)\le f_3(\xi)=\frac{2(C+1)^2|\xi_0|^2}{|\xi_0|^2+(|\xi_1|^2+|\xi_3|^2)^2}\in A_{X_0}.
\end{equation}
Taking into account monotonicity of the Gelfand transform (cf. Lemma~\ref{Lem2}), we derive from (\ref{15}) that
$$
0<|\tilde\xi_{X_1}|^2|\tilde\xi_{X_2}|^2=\widehat{f_1}(\eta)\widehat{f_2}(\eta)\le\widehat{f_3}(\eta)=
2(C+1)^2|\tilde\xi_{X_0}|^2.
$$
Hence, $\tilde\xi_{X_0}\not=0$ and $X_0=X_1\cap X_2\in\tilde\cL$.
Let $\tilde X$ be a subspace from $\tilde\cL$ of minimal dimension. As was already established, for each $X\in\tilde\cL$ the subspace $X_0=X\cap\tilde X\in\tilde\cL$. Since $X_0\subset\tilde X$ while $\dim\tilde X\le\dim X_0$, we obtain that $\tilde X=X_0\subset X$. Thus, $X\supset\tilde X$ $\forall X\in\tilde\cL$. Let us demonstrate that, conversely, any subspace $X\supset\tilde X$ belongs to $\tilde L$ and $\tilde\xi_X\upuparrows\tilde\xi_{\tilde X}$. For that, we introduce the space $X_1=X\ominus\tilde X$,
so that $\R^n=\tilde X\oplus X_1\oplus X^\bot$. Denote by $\xi_0$, $\xi_1$, $\xi_2$ the orthogonal projections of a vector $\xi\in\R^n$ on the subspaces $\tilde X$, $X_1$, $X^\bot$, respectively. Then $\xi=\xi_0+\xi_1+\xi_2$.
For arbitrary $u,v\in\R^n$ we find
\begin{eqnarray*}
\frac{u\cdot\xi_0}{(|\xi_0|^2+(|\xi_1|^2+|\xi_2|^2)^2)^{1/2}}\cdot
\frac{v\cdot(\xi_0+\xi_1)}{(|\xi_0|^2+|\xi_1|^2+|\xi_2|^4)^{1/2}}=\\
\frac{v\cdot\xi_0}{(|\xi_0|^2+(|\xi_1|^2+|\xi_2|^2)^2)^{1/2}}\cdot
\frac{u\cdot\xi_0}{(|\xi_0|^2+|\xi_1|^2+|\xi_2|^4)^{1/2}}+\gamma(\xi),
\end{eqnarray*}
where $\displaystyle\gamma(\xi)=\frac{v\cdot\xi_1}{(|\xi_0|^2+(|\xi_1|^2+|\xi_2|^2)^2)^{1/2}}\in A_0$.
Applying the Gelfand transform to the above equality, we obtain the equality
\begin{equation}\label{16}
(u\cdot\tilde\xi_{\tilde X})(v\cdot\tilde\xi_X)=(u\cdot\pr_{\tilde X}\tilde\xi_X)(v\cdot\tilde\xi_{\tilde X}).
\end{equation}
Taking $u=\tilde\xi_{\tilde X}$, $v\perp\tilde X$, we derive from (\ref{16}) that $v\cdot\tilde\xi_X=0$
for all $v\perp\tilde X$, which implies the inclusion $\tilde\xi_X\in\tilde X$. In particular, $\pr_{\tilde X}\tilde\xi_X=\tilde\xi_X$ and it follows from (\ref{16}) that
$$
(u\cdot\tilde\xi_{\tilde X})(v\cdot\tilde\xi_X)=(u\cdot\tilde\xi_X)(v\cdot\tilde\xi_{\tilde X}) \quad \forall u,v\in\R^n.
$$
In view of this relation we find that $\tilde\xi_X=c\tilde\xi_{\tilde X}$ for some real constant $c$.
Further,
\begin{equation}\label{17}
\frac{|\xi_0|^2}{|\xi_0|^2+(|\xi_1|^2+|\xi_2|^2)^2}=\frac{|\xi_0|^2}{|\xi_0|^2+|\xi_1|^2+|\xi_2|^4}\cdot
\frac{|\xi_0|^2+|\xi_1|^2+|\xi_2|^4}{|\xi_0|^2+(|\xi_1|^2+|\xi_2|^2)^2}.
\end{equation}
Observe that
$$
\frac{|\xi_0|^2+|\xi_1|^2+|\xi_2|^4}{|\xi_0|^2+(|\xi_1|^2+|\xi_2|^2)^2}=
g(\xi)=\frac{|\xi_0|^2+|\xi_2|^4}{|\xi_0|^2+(|\xi_1|^2+|\xi_2|^2)^2}
$$
up to a term vanishing at infinity, and $g(\xi)\in A_{\tilde X}$. Hence, applying the Gelfand transform to (\ref{17}), we obtain
$$
0<|\tilde\xi_{\tilde X}|^2=|\tilde\xi_X|^2\hat g(\eta).
$$
It follows from this relation that $\tilde\xi_X\not=0$, and the constant $c\not=0$. Finally,
$c|\tilde\xi_{\tilde X}|^2=\tilde\xi_X\cdot\tilde\xi_{\tilde X}=\hat h(\eta)$, where
$$
h(\xi)=\frac{|\xi_0|^2}{(|\xi_0|^2+|\xi_1|^2+|\xi_2|^4)^{1/2}(|\xi_0|^2+(|\xi_1|^2+|\xi_2|^2)^2)^{1/2}}\ge 0.
$$
By the monotonicity of the  Gelfand transform, we find that $c>0$. Therefore, $\tilde\xi_X\upuparrows\tilde\xi_{\tilde X}$. Denote $\zeta_1=\tilde\xi_{\tilde X}/|\tilde\xi_{\tilde X}|\in\tilde X$ (~remark that $\zeta_1=\tilde\xi_{\R^n}$~). Thus,
\begin{equation}\label{18}
\tilde\xi_X\not=0\Leftrightarrow X\supset \tilde X  \ \mbox{ and } \tilde\xi_X\upuparrows\zeta_1.
\end{equation}

{\bf 2nd Step.}

We introduce the family $\bar\cL=\{ \ X\subset\R^n \ | \ \bar\xi_X=0 \}$. Let $X_1,X_2\in\bar\cL$. We show that $X_0=X_1\cap X_2\in\bar\cL$. For that,  we denote $X_{10}=X_1\ominus X_0$, $X_{20}=X_2\ominus X_0$. Then representations (\ref{10}) and (\ref{11})
hold. We introduce the functions
$$
g_1(\xi)=\frac{|\xi_3|^4}{|\xi_0|^2+|\xi_1|^2+|\xi_3|^4}, \quad g_2(\xi)=\frac{|\xi_4|^4}{|\xi_0|^2+|\xi_2|^2+|\xi_4|^4},
$$
and remark that $\widehat{g_1}(\eta)=\widehat{g_2}(\eta)=0$, in view of the condition $\bar\xi_{X_1}=\bar\xi_{X_2}=0$.
Since
\begin{eqnarray*}
h_1(\xi)=\frac{|\xi_3|^4}{|\xi_0|^2+(|\xi_1|^2+|\xi_3|^2)^2}=g_1(\xi)\frac{|\xi_0|^2+|\xi_1|^2+|\xi_3|^4}{|\xi_0|^2+(|\xi_1|^2+|\xi_3|^2)^2},\\
h_2(\xi)=\frac{|\xi_4|^4}{|\xi_0|^2+(|\xi_2|^2+|\xi_4|^2)^2}=g_2(\xi)\frac{|\xi_0|^2+|\xi_2|^2+|\xi_4|^4}{|\xi_0|^2+(|\xi_2|^2+|\xi_4|^2)^2},
\end{eqnarray*}
then $0\le h_k(\xi)\le 2g_k(\xi)$ for $k=1,2$ and sufficiently large $|\xi|$. Therefore, $0\le \widehat{h_k}(\eta)\le 2\widehat{g_k}(\eta)=0$, $k=1,2$, and we arrive at $\widehat{h_k}(\eta)=0$ for $k=1,2$. Remark also that $|\xi_1|^2+|\xi_3|^2=|\xi_2|^2+|\xi_4|^2=|\xi-\xi_0|^2$, and, therefore,
\begin{equation}\label{18a}
p(\xi)\doteq |\xi_0|^2+(|\xi_1|^2+|\xi_3|^2)^2=|\xi_0|^2+(|\xi_2|^2+|\xi_4|^2)^2.
\end{equation}
By estimates (\ref{est}) we see that $|\xi_1|^2\le C(|\xi_3|^2+|\xi_4|^2)$. This inequality together with (\ref{18a}) imply that
$$
\frac{|\xi_1|^2}{(p(\xi))^{1/2}}\le C\frac{|\xi_3|^2}{(p(\xi))^{1/2}}+C\frac{|\xi_4|^2}{(p(\xi))^{1/2}}=C(h_1(\xi))^{1/2}+C(h_2(\xi))^{1/2}.
$$
Therefore,
$$
h(\xi)\doteq\frac{|\xi_1|^2+|\xi_3|^2}{(p(\xi))^{1/2}}\le (C+1)(h_1(\xi))^{1/2}+C(h_2(\xi))^{1/2}.
$$
This implies that $|\bar\xi_{X_0}|^2=\hat h(\eta)\le (C+1)(\widehat{h_1}(\eta))^{1/2}+C(\widehat{h_2}(\eta))^{1/2}=0$.
Hence, $\bar\xi_{X_0}=0$ and $X_0\in\bar\cL$. This statement allows to establish existence of minimal element $\bar X$ in $\bar\cL$,
in the same way as for the family $\tilde\cL$. Namely, let $\bar X$ be an element in $\bar\cL$ of minimal dimension. Then for arbitrary
$X\in\bar\cL$ the intersection $X_0=\bar X\cap X\in\bar\cL$. Since $X_0\subset\bar X$ while $\dim\bar X\le\dim X_0$, then $\bar X=X_0\subset X$. Hence $\bar X$ is the smallest subspace in $\bar\cL$. Notice also that if a subspace $X\supset\bar X$ then
$X\in\bar\cL$. Indeed, $\R^n=\bar X\oplus X_1\oplus X^\bot$, where $X_1=X\ominus\bar X$. Therefore, $\xi=\xi_0+\xi_1+\xi_2$,
with $\xi_0,\xi_1,\xi_2$ being the orthogonal projection of $\xi\in\R^n$ on the subspaces $\bar X, X_1, X^\bot$, respectively.
Let
$$
\rho(\xi)=\frac{|\xi_0|^2+|\xi_1|^2+|\xi_2|^4}{|\xi_0|^2+(|\xi_1|^2+|\xi_2|^2)^2}.
$$
Then
$$
\hat\rho(\eta)=|\tilde\xi_{\bar X}|^2+|\pr_{X^\bot}\bar\xi_{\bar X}|^4=|\tilde\xi_{\bar X}|^2=1
$$
because $\bar\xi_{\bar X}=0$ while $|\tilde\xi_{\bar X}|^2+|\bar\xi_{\bar X}|^4=1$.

Let $\displaystyle q(\xi)=\frac{|\xi_2|^4}{|\xi_0|^2+|\xi_1|^2+|\xi_2|^4}$. Since
$\displaystyle\rho(\xi)q(\xi)=\frac{|\xi_2|^4}{|\xi_0|^2+(|\xi_1|^2+|\xi_2|^2)^2}$, we find
$\hat q(\eta)=\hat\rho(\eta)\hat q(\eta)=|\pr_{X^\bot}\bar\xi_{\bar X}|^4=0$, which implies that $\bar\xi_X=0$.
Thus, $X\in\bar\cL$, and
\begin{equation}\label{19}
\bar\cL=\{ \ X\subset\R^n \ | \ X\supset\bar X \ \}.
\end{equation}
Notice that $|\tilde\xi_{\bar X}|=1$. Therefore, $\bar X\in\tilde\cL$ and, in view of (\ref{18}), $\bar X\supset\tilde X$.

{\bf 3rd Step.}

Assume that $X_1\subset X_2\subset\R^n$ and $\bar\xi_{X_2}\not=0$. We claim that $\bar\xi_{X_1}\not=0$ and
$\bar\xi_{X_2}\upuparrows\zeta\doteq\pr_{X_2^\bot}\bar\xi_{X_1}$ (that is, $\zeta=c\bar\xi_{X_2}$ for some $c\ge 0$).

Indeed, if $\bar\xi_{X_1}=0$ then $X_1\in\bar\cL$. By (\ref{19}) we find
$X_2\in\bar\cL$. But this contradicts to the assumption $\bar\xi_{X_2}\not=0$. Further, let
$$
p_1(\xi)=(|\xi_1|^2+(|\xi_2|^2+|\xi_3|^2)^2)^{1/4}, \quad
p_2(\xi)=(|\xi_1|^2+|\xi_2|^2+|\xi_3|^4)^{1/4},
$$
where $\xi_1=\pr_{X_1}\xi$, $\xi_2=\pr_{X_2\ominus X_1}\xi$, $\xi_3=\pr_{X_2^\bot}\xi$.
Evidently, for each $u,v\in\R^n$
$$
\frac{u\cdot\xi_3}{p_1(\xi)}\cdot\frac{v\cdot\xi_3}{p_2(\xi)}=\frac{v\cdot\xi_3}{p_1(\xi)}\cdot\frac{u\cdot\xi_3}{p_2(\xi)}.
$$
Applying the Gelfand transform to this identity, we find
\begin{equation}\label{20}
(u\cdot\zeta)(v\cdot\bar\xi_{X_2})=(v\cdot\zeta)(u\cdot\bar\xi_{X_2}) \ \forall u,v\in\R^n,
\end{equation}
where $\zeta=\pr_{X_2^\bot}\bar\xi_{X_1}$. It readily follows from (\ref{20}) that $\zeta=c\bar\xi_{X_2}$ for some constant $c\in\R$.
Since $\zeta\cdot\bar\xi_{X_2}$ coincides with the Gelfand transform of the nonnegative symbol
$\displaystyle\frac{|\xi_3|^2}{p_1(\xi)p_2(\xi)}$, we conclude that $\zeta\cdot\bar\xi_{X_2}\ge 0$, i.e. $c\ge 0$.

{\bf 4th Step.}

In this step we prove that for any $X\subset\bar X$ the vector $\bar\xi_X\in\bar X$. Moreover, in the case $X\subset\tilde X$, $X\not=\tilde X$ the vector $\bar\xi_X\in\tilde X$.

First, we notice that if $X=\bar X$, then $\bar\xi_X=0\in\bar X$. In the remaining case $X\not=\bar X$ we denote by $\xi_1,\xi_2,\xi_3$ the orthogonal projections of $\xi\in\R^n$ onto the subspaces $X$, $\bar X\ominus X$, $\bar X^\bot$, respectively, and introduce the symbols
$$
a(\xi)=\frac{|\xi_3|^4}{|\xi_1|^2+(|\xi_2|^2+|\xi_3|^2)^2}\in A_X, \
b(\xi)=\frac{|\xi_3|^4}{|\xi_1|^2+|\xi_2|^2+|\xi_3|^4}\in A_{\bar X}.
$$
as is easy to verify, $a(\xi)\le 2b(\xi)$ for sufficiently large $|\xi|$, which implies
$$
|\pr_{\bar X^\bot}\bar\xi_X|^4=a(\eta)\le 2\hat b(\eta)=2|\bar\xi_{\bar X}|^4=0 \Rightarrow \pr_{\bar X^\bot}\bar\xi_X=0.
$$
This means that $\bar\xi_X\in\bar X$, as was to be proved.

It remains only to consider the case when $X\varsubsetneq\tilde X$. Let $\xi_1,\xi_2,\xi_3$ be the orthogonal projections of $\xi\in\R^n$ onto the subspaces $X$, $\tilde X\ominus X$, $\tilde X^\bot$, respectively. We introduce the functions
\begin{eqnarray*}
a(\xi)=\frac{|\xi_1|^2}{|\xi_1|^2+(|\xi_2|^2+|\xi_3|^2)^2}\in A_X, \
b(\xi)=\frac{|\xi_1|^2}{|\xi_1|^2+|\xi_2|^2+|\xi_3|^4}\in A_{\tilde X}, \\
c(\xi)=\frac{|\xi_1|^2+|\xi_2|^2+|\xi_3|^4}{|\xi_1|^2+(|\xi_2|^2+|\xi_3|^2)^2}\sim
\frac{|\xi_1|^2+|\xi_3|^4}{|\xi_1|^2+(|\xi_2|^2+|\xi_3|^2)^2}.
\end{eqnarray*}
Since $X\notin\tilde\cL$, then $\tilde\xi_X=0$ and
\begin{equation}\label{21}
\hat c(\eta)=|\tilde\xi_X|^2+|\pr_{\tilde X^\bot}\bar\xi_X|^4=|\pr_{\tilde X^\bot}\bar\xi_X|^4.
\end{equation}
Further, $\hat a(\eta)=|\tilde\xi_X|^2=0$, $\hat b(\eta)=|\tilde\xi_{\tilde X}|^2\not=0$, and $0=\hat a(\eta)=\hat b(\eta)\hat c(\eta)$. Therefore, $\hat c(\eta)=0$, and it follows from (\ref{21}) that $\pr_{\tilde X^\bot}\bar\xi_X=0$, that is,
$\bar\xi_X\in\tilde X$.

{\bf 5th Step.}

Here we construct the orthonormal family $\{\zeta_k\}_{k=1}^m$.
First, we set $\zeta_1=\tilde\xi_{\R^n}=\bar\xi_{\{0\}}$. Assuming that the vectors $\zeta_1,\ldots,\zeta_{k-1}$ have already known, we define \begin{equation}\label{22}
\zeta_k=\bar\xi_{X_{k-1}}/|\bar\xi_{X_{k-1}}|\in X_{k-1}^\bot,
\end{equation}
 where $X_{k-1}$ is a subspace spanned by the vectors $\zeta_1,\ldots,\zeta_{k-1}$ (notice that $X_0=\{0\}$). This definition is correct while $X_{k-1}\varsubsetneq\bar X$
because by (\ref{19}) $\bar\xi_{X_{k-1}}\not=0$. As was demonstrated in the 4th step, $\zeta_k\in\bar X$.
We see that the construction of $\zeta_k$ may be continued until $k=m=\dim\bar X$. The $m$-dimensional subspace $X_m\subset\bar X$ must coincide with $\bar X$: $X_m=\bar X$, so that $\bar\xi_{X_m}=0$. By the construction
$\{\zeta_k\}_{k=1}^m$ is an orthonormal basis in $\bar X$. Let $d=\dim\tilde X$. Then $1\le d\le m$. By the second statement proven in 4th Step $\zeta_k\in\tilde X$ while $X_{k-1}\varsubsetneq\tilde X$. Since $\zeta_1\in\tilde X$
then by induction $X_k\subset \tilde X$ for all $1\le k\le d$. Comparing the dimension, we claim that $\tilde X=X_d=
\cL(\zeta_1,\ldots,\zeta_d)$. As was shown in 1st step, for $\tilde\xi_X\not=0$ this vector is co-directed with $\zeta_1$. The proof of (i) is complete.

To complete the proof of statement (ii), we choose a subspace $X\subset\R^n$ such that $\bar\xi_X\not=0$. Then,
in view of (\ref{19}) $X\not\supset\bar X=\cL(\zeta_1,\ldots,\zeta_m)$, Therefore, there exists the vector
$\zeta_k\notin X$. Let $k=k(X)=\min\{ \ k=1,\ldots,m \ | \ \zeta_k\notin X \ \}$. Then $X_{k-1}\subset X$,
$\zeta_k\notin X$. Since $\zeta_k\upuparrows\bar\xi_{X_{k-1}}\not=0$, then by the assertion established in
the 3rd Step we claim that $\bar\xi_X\upuparrows\pr_{X^\bot}\zeta_k\not=0$, as was to be proved.

Remark also that by results of the 3rd Step requirement (ii) for $X=X_{k-1}$ implies (\ref{22}). This readily implies that
the orthonormal family $\zeta_k$, $k=1,\ldots,m$ is uniquely defined by the point $\eta$. The parameter $d$ is also uniquely determined by the condition $d=\dim\tilde X$.

{\bf 6th Step.} It only remains to show that $d\ge m-1$. Assuming the contrary $d\le m-2$, we see that the space $X_1$ spanned by the vectors $\xi_k$, $k=1,\ldots,d+1$ is a proper subspace of $\bar X$: $\tilde X\varsubsetneq X_1\varsubsetneq\bar X$. We extend the system $\zeta_k$, $k=1,\dots,m$ to an orthonormal basis $\zeta_k$, $k=1,\ldots,n$ in $\R^n$. Let $s_k=s_k(\xi)$, $k=1,\ldots,n$ be coordinates of a vector $\xi\in\R^n$ in this basis: $\displaystyle\xi=\sum_{k=1}^n s_k\zeta_k$. We introduce the following functions
\begin{eqnarray*}
p_1(\xi)=\frac{s_1^2}{\sum_{k=1}^d s_k^2+\left(\sum_{k=d+1}^n s_k^2\right)^2}, \
q_1(\xi)=\frac{s_{d+2}^4}{\sum_{k=1}^{d+1} s_k^2+\left(\sum_{k=d+2}^n s_k^2\right)^2}, \\
p_2(\xi)=\frac{s_1^2}{\sum_{k=1}^{d+1} s_k^2+\left(\sum_{k=d+2}^n s_k^2\right)^2}, \
q_2(\xi)=\frac{s_{d+2}^4}{\sum_{k=1}^d s_k^2+\left(\sum_{k=d+1}^n s_k^2\right)^2}.
\end{eqnarray*}
Obviously, $p_1,q_2\in S_{\tilde X}$, $p_2,q_1\in S_{X_1}$, and $p_1q_1=p_2q_2$. Therefore,
\begin{equation}\label{23}
\widehat{p_1}(\eta)\widehat{q_1}(\eta)=\widehat{p_2}(\eta)\widehat{q_2}(\eta).
\end{equation}
Now observe that $\widehat{p_1}(\eta)=|\tilde \xi_{\tilde X}|^2\not=0$, $\widehat{p_2}(\eta)=|\tilde\xi_{X_1}|^2\not=0$,
$\widehat{q_1}(\eta)=|\bar\xi_{X_1}|^4\not=0$ (because $X_1\varsubsetneq\bar X$), and
$\widehat{q_2}(\eta)=|\bar\xi_{\tilde X}\cdot\zeta_{d+2}|^4=0$ because $\bar\xi_{\tilde X}\upuparrows\xi_{d+1}\perp\xi_{d+2}$. Hence $\widehat{p_1}(\eta)\widehat{q_1}(\eta)\not=0$,
$\widehat{p_2}(\eta)\widehat{q_2}(\eta)=0$, which contradicts (\ref{23}). The proof is complete.
\end{pf}

The statement of Proposition~\ref{Pro3} is sharp, in the sense that for every orthonormal system $\zeta_k$, $k=1,\ldots,m$, and an integer number $d\in\{m-1,m\}$, one can find a point $\eta\in\S$ such that the statements (i), (ii) of Proposition~\ref{Pro4} hold. To prove this assertion, we need the notion of an essential ultrafilter.
We call sets $A,B\subset\R^n$ equivalent: $A\sim B$ if $\mu(A\vartriangle B)=0$, where $A\vartriangle B=(A\setminus B)\cup (B\setminus A)$ is the symmetric difference and $\mu$ is the outer Lebesgue measure.
Let $\mathfrak{F}$ be a filter in $\R^n$. This filter is called \textit{essential} if from the conditions
$A\in\mathfrak{F}$ and $B\sim A$ it follows that $B\in A$. It is clear that an essential filter cannot include sets
of null measure, since such sets are equivalent to $\emptyset$.
Using Zorn's lemma, one can prove that any essential filter is contained in a maximal essential filter. Maximal
essential filters are called essential ultrafilters.

\begin{lem}\label{Lem3} Let $\mathfrak{U}$ be an essential ultrafilter. Then for each $A\subset\R^n$ either $A\in\mathfrak{U}$ or $\R^n\setminus A\in\mathfrak{U}$.
\end{lem}

\begin{pf} Assuming that $A\notin\mathfrak{U}$, we introduce
$$\mathfrak{F}=\{ \ B\subset\R^n \ | \ B\cup A\in\mathfrak{U} \ \}.$$
Obviously, $\mathfrak{F}$ is an essential filter, $\R^n\setminus A\in\mathfrak{F}$, and $\mathfrak{U}\le\mathfrak{F}$. Since the filter
$\mathfrak{U}$ is maximal, we obtain that $\mathfrak{U}=\mathfrak{F}$. Hence, $\R^n\setminus A\in\mathfrak{U}$. The proof is complete.
\end{pf}

The property indicated in Lemma~\ref{Lem3} is the characteristic property of ultrafilters, see for example, \cite{Burb}. Therefore, we have the following statement.

\begin{cor}\label{Cor1} Any essential ultrafilter is an ultrafilter, i.e. a maximal element in a set of all filters. \end{cor}

\begin{lem}\label{Lem4} Let $\mathfrak{U}$ be an essential ultrafilter, and $f(\xi)$ be a bounded function in $\R^n$. Then there exists $\displaystyle\lim_{\mathfrak{U}} f(\xi)$. If a function $g(\xi)=f(\xi)$ almost everywhere on $\R^n$, then there exists $\displaystyle\lim_{\mathfrak{U}} g(\xi)=\lim_{\mathfrak{U}} f(\xi)$.
\end{lem}
\begin{pf}
By Corollary~\ref{Cor1} \ $\mathfrak{U}$ is an ultrafilter. By the known properties of ultrafilters, the image $f_*\mathfrak{U}$ is an ultrafilter on the compact $[-M,M]$, where $M=\sup|f(\xi)|$, and this ultrafilter converges to some point $x\in [-M,M]$. Therefore, $\displaystyle\lim_{\mathfrak{U}} f(\xi)=\lim f_*\mathfrak{U}=x$. Further, suppose that a function $g=f$ a.e. on $\R^n$. Then the set $E=\{ \xi\in\R^n \ | \ g(\xi)\not=f(\xi) \ \}$ has null Lebesgue measure. Let $V$ be a neighborhood of $x$. Then $g^{-1}(V)\supset f^{-1}(V)\setminus E$. By the convergence of the ultrafilter $f_*\mathfrak{U}$ the set $f^{-1}(V)\in\mathfrak{U}$. Since $\mathfrak{U}$ is an essential ultrafilter while $f^{-1}(V)\setminus E\sim f^{-1}(V)$, then $f^{-1}(V)\setminus E\in\mathfrak{U}$. This set is contained in $g^{-1}(V)$, and we claim that $g^{-1}(V)\in\mathfrak{U}$. Since $V$ is an arbitrary neighborhood of $x$, we conclude that $\displaystyle\lim_{\mathfrak{U}} g(\xi)=x$. The proof is complete.
\end{pf}

By the statement of Lemma~\ref{Lem4}, the functional $\displaystyle f\to\lim_{\mathfrak{U}} f(\xi)$ is well-defined on $L^\infty(\R^n)$ and it is a linear multiplicative functional on $L^\infty(\R^n)$. In other words, this functional belongs to the spectrum of algebra $L^\infty(\R^n)$ (actually, this spectrum coincides with the space of such functionals).

Now we are ready to prove the sharpness of Proposition~\ref{Pro3}.

\begin{prop}\label{Pro4} Let $\zeta_k$, $k=1,\ldots,m$ be an orthonormal system in $\R^n$, $1\le d\in\{m-1,m\}$. Then there exists a point $\eta\in\S$ such that the statements (i), (ii) of Proposition~\ref{Pro3} hold.
\end{prop}
\begin{pf}
We extend vectors $\zeta_k$, $k=1,\ldots,m$ to a basis $\zeta_k$, $k=1,\ldots,n$ in $\R^n$. Let
$\sigma_k$, $k=2,\ldots,n$ be a decreasing family of positive numbers such that $1>\sigma_2>\cdots>\sigma_d>1/2\ge\sigma_{d+1}>\cdots>\sigma_n>0$, and $\sigma_{d+1}=1/2$ only if $d=m-1$.
We introduce the sets
$$B_r=\left\{ \ \xi=\sum_{k=1}^n s_k\zeta_k \ | \ |\xi|>r, s_1>0, \ \mbox{ and } s_1^{\sigma_k}<s_k<2s_1^{\sigma_k} \ \forall k=2,\ldots, n \ \right\},$$
$r>0$. It is clear that $B_r$ are nonempty open sets in $\R^n$, which form the base of some essential filter $\mathfrak{F}$. Let $\mathfrak{U}$ be an essential ultrafilter such that $\mathfrak{F}\le\mathfrak{U}$. Since the limit along $\mathfrak{U}$ is a linear multiplicative functional on $A$ vanishing on the ideal $A_0$, it forms a linear multiplicative functional on $\A$, and there exists a unique element $\eta\in\S$ such that
$\displaystyle\hat a(\eta)=\lim_{\mathfrak{U}}a(\xi)$ for each $a\in\A$. We will demonstrate that the element $\eta$ satisfies conditions
(i), (ii) of Proposition~\ref{Pro4}. Assume that a subspace $X\not\supset\tilde X=\cL(\zeta_1,\ldots,\zeta_d)$. Then there exists $k$, $1\le k\le d$ such that $\zeta_k\notin X$. Let $k=k(X)$ be the minimal one among such $k$.
We denote by $P_1$, $P_2$ the orthogonal projections onto the spaces $X$, $X^\bot$, respectively, and set $v_i=P_2\zeta_i$. Then $v_k\not=0$ while $v_i=0$ for $1\le i<k$. If $\xi\in B_r$, then
\begin{equation}\label{24}
r^2<|\xi|^2=\sum_{i=1}^n s_i^2\le s_1^2+2\sum_{i=2}^n s_1^{2\sigma_i}\le C_rs_1^2,
\end{equation}
where $C_r\to 1$ as $r\to\infty$. Here we take into account the condition $\sigma_i<1$. In particular, it follows from (\ref{24}) that $s_1>r/2$ for large $r$. Denote, as above, $\tilde\xi=P_1\xi$, $\bar\xi=P_2\xi$. Since $\sigma_i<\sigma_k$ for $i>k$, $s_1>r/2\mathop{\to}\limits_{r\to\infty}\infty$, and $|v_k|>0$, we find that for sufficiently large $r$
\begin{equation}\label{25}
|\bar\xi|=\left|\sum_{i=k}^n s_iv_i\right|\ge s_k|v_k|-\sum_{i=k+1}^ns_i|v_i|\ge s_1^{\sigma_k}|v_k|-
2\sum_{i=k+1}^ns_1^{\sigma_i}|v_i|\ge cs_1^{\sigma_k},
\end{equation}
where  $c=\const>0$. It follows from (\ref{24}), (\ref{25}) that for $\xi\in B_r$, where $r$ is sufficiently large
\begin{equation}\label{26}
a(\xi)=\frac{|\tilde\xi|^2}{|\tilde\xi|^2+|\bar\xi|^4}\le \frac{|\xi|^2}{|\bar\xi|^4}\le\frac{C_r}{c^4}s_1^{2-4\sigma_k}\to 0
\end{equation}
as $r\to\infty$ because $\sigma_k>1/2$ and $s_1\to \infty$ as $r\to\infty$.
It follows from (\ref{26}) that
$$
|\tilde\xi_X(\eta)|^2=\hat a(\eta)=\lim_{\mathfrak{U}}a(\xi)=\lim_{\mathfrak{F}}a(\xi)=0.
$$
We claim that $\tilde\xi_X(\eta)=0$.

If $X\supset\tilde X$, then $\displaystyle\bar\xi=\sum_{i=d+1}^n s_iv_i$, where $v_i=P_2\zeta_i$, and for $\xi\in B_r$
\begin{eqnarray}\label{27}
s_1^2\le |\tilde\xi|^2\le |\xi|^2=\sum_{i=1}^ns_i^2\le C_1s_1^2, \\
\label{28}
|\bar\xi|\le  \sum_{i=d+1}^n s_i\le C_2s_1^{\sigma_{d+1}}, \quad C_1,C_2=\const.
\end{eqnarray}
Since $\sigma_{d+1}\le 1/2$, then it follows from (\ref{27}), (\ref{28}) that for sufficiently large $r$
$$
a(\xi)=\frac{|\tilde\xi|^2}{|\tilde\xi|^2+|\bar\xi|^4}\ge (C_1+C_2^4)^{-1}>0,
$$
which implies that
$|\tilde\xi_X(\eta)|^2=\hat a(\eta)=\lim\limits_{\mathfrak{U}}a(\xi)>0$. Hence $\tilde\xi_X=\tilde\xi_X(\eta)\not=0$.
Observe also that, as follows from (\ref{27}), (\ref{28}),
$$
\frac{s_i}{(|\tilde\xi|^2+|\bar\xi|^4)^{1/2}}\le cs_1^{\sigma_2-1}\mathop{\to}_{r\to\infty} 0, \quad
i=2,\ldots,n,
$$
and, therefore,
$$\tilde\xi_X=\lim_{\mathfrak{U}}\frac{\tilde\xi}{(|\tilde\xi|^2+|\bar\xi|^4)^{1/2}}=
\left(\lim_{\mathfrak{U}}\frac{s_1}{(|\tilde\xi|^2+|\bar\xi|^4)^{1/2}}\right)\zeta_1\upuparrows\zeta_1.
$$
We conclude that condition (i) is satisfied.

To prove (ii), assume that $X\supset\bar X=\cL(\zeta_1,\ldots,\zeta_m)$. Then
$\displaystyle\bar\xi=\sum_{i=m+1}^n s_iv_i$, $v_i=P_2\zeta_i$, which implies the estimate
$$
|\bar\xi|\le 2\sum_{i=m+1}^ns_1^{\sigma_i}\le Cs_1^{\sigma_{m+1}}, \ C=\const,
$$
for all $\xi\in B_r$ with sufficiently large $r$. On the other hand,
$|\tilde\xi|\ge s_1$. Therefore, for $\xi\in B_r$
$$
\frac{|\bar\xi|^4}{|\tilde\xi|^2+|\bar\xi|^4}\le C^4s_1^{4\sigma_{m+1}-2}\mathop{\to}_{r\to\infty} 0
$$
because $\sigma_{m+1}<1/2$. Hence,
$$
|\bar\xi_X|^4=\lim_{\mathfrak{U}}\frac{|\bar\xi|^4}{|\tilde\xi|^2+|\bar\xi|^4}=0,
$$
that is, $\bar\xi_X=0$.

Now, suppose that $X\not\supset\bar X$. Then there exists $\zeta_k\notin X$, where $1\le k\le m$. We chose
$k$ being the minimal one. Then $\zeta_i\in X$, $1\le i<k$, and $\displaystyle\bar\xi=\sum_{i=k}^n s_iv_i$, $v_i=P_2\zeta_i$, which implies the estimate
\begin{equation}\label{29}
|\bar\xi|\ge s_k|v_k|-\sum_{i=k+1}^ns_i|v_i|\ge s_1^{\sigma_k}|v_k|-2\sum_{i=k+1}^ns_1^{\sigma_i}|v_i|\ge cs_1^{\sigma_k}, \ c=|v_k|/2>0,
\end{equation}
for all $\xi\in B_r$ with sufficiently large $r$. We use here that $v_k\not=0$ and $\sigma_k>\sigma_i$ for $i>k$. Further,
\begin{equation}\label{30}
|\tilde\xi|^2\le|\xi|^2=\sum_{i=1}^ns_i^2\le s_1^2+2\sum_{i=2}^ns_1^{2\sigma_i}\le 2|s_1|^2
\end{equation}
for all $\xi\in B_r$ with large $r$. It follows from (\ref{29}), (\ref{30}) and from the condition $\sigma_k\ge\sigma_m\ge 1/2$ that
\begin{equation}\label{31}
cs_1^{\sigma_k}\le(|\tilde\xi|^2+|\bar\xi|^4)^{1/4}\le Cs_1^{\sigma_k}, \quad C=\const.
\end{equation}
In view of (\ref{31}) for all $\xi\in B_r$ with sufficiently large $r$
\begin{eqnarray*}
\frac{s_i}{(|\tilde\xi|^2+|\bar\xi|^4)^{1/4}}\le \frac{2}{c}s_1^{\sigma_i-\sigma_k}\mathop{\to}_{r\to\infty}0, \quad k+1\le i\le n, \\
\frac{1}{C}\le\frac{s_k}{(|\tilde\xi|^2+|\bar\xi|^4)^{1/4}}\le \frac{2}{c}.
\end{eqnarray*}
This implies that
$$
\bar\xi_X(\eta)=av_k\upuparrows \pr_{X^\bot}\zeta_k,
$$
where
$$
a=\lim_{\mathfrak{U}}\frac{s_k}{(|\tilde\xi|^2+|\bar\xi|^4)^{1/4}}>0,
$$
and
$k=k(X)=\min\{ \ k=1,\ldots,m \ | \ \zeta_k\notin X \ \}$.

We see that requirement (ii) of Proposition~\ref{Pro4} is also satisfied.
The proof is complete.
\end{pf}

\begin{rem}\label{Rem1} For each real $t\ne 0$ the map $h_t(\psi)(\xi)=\psi^t\doteq\psi(t\xi)$ is an isomorphism of algebra $\A$. Indeed, it is easy to verify that
$$
A_{\psi^t}=Q_t^{-1}A_\psi Q_t, \quad B_{\Phi}=Q_t^{-1}B_{\Phi^t} Q_t \quad \forall\psi(\xi)\in A, \Phi(x)\in C_0(\R^n),
$$
where $\Phi^t(x)=\Phi(tx)$, and the operator $Q_t$ in $L^2(\R^n)$ is defined by the equality $Q_tu(x)=u^t=u(tx)$.
Therefore, the operators
$[A_{\psi^t},B_{\Phi}]=Q_t^{-1}[A_\psi,B_{\Phi^t}]Q_t$ are compact in $L^2$ for all $\Phi(x)\in C_0(\R^n)$. This implies that $h_t$ is well-defined on $A$ and evidently transfers the ideal $A_0$ into itself. This allows to define the operator $h_t$ on the quotient algebra $\A=A/A_0$. It is clear that $h_t$ is invertible and $h_t^{-1}=h_{1/t}$.
Therefore, the operator $h_t$ generates the corresponding homeomorphism of the spectrum $\widehat{h_t}:\S\to\S$, so that
$\widehat{\psi}(\widehat{h_t}(\eta))=\widehat{h_t(\psi)}(\eta)$. We denote $\widehat{h_t}(\eta)=t\eta$.
This determines an action of the multiplicative group of $\R$ on the space $\S$. If $X$ is a subspace of $\R^n$,
and $(\tilde\xi(\eta),\bar\xi(\eta))=p_X(\eta)$, then it is directly verified that for each $t\ne 0$
$$
\tilde\xi(t\eta)=a(t,\eta)\tilde\xi(t\eta), \ \bar\xi(t\eta)=b(t,\eta)\bar\xi(t\eta),
$$
where
$$
a(t,\eta)=t(t^2|\tilde\xi(\eta)|^2+t^4|\bar\xi(\eta)|^4)^{-1/2}, \ b(t,\eta)=t(t^2|\tilde\xi(\eta)|^2+t^4|\bar\xi(\eta)|^4)^{-1/4}.
$$
In particular, $(b(t,\eta))^2=ta(t,\eta)$.
\end{rem}

\section{$H$-measures and the localization property}

Now, let $\Omega\subset\R^n$ be an open domain and $U_r(x)\in L^2_{loc}(\Omega,\C^N)$ be a sequence of generally complex-valued vector functions weakly convergent to the zero vector. Denote by
$\Blim$ a generalized Banach limit (see \cite{Ban}), that is, a linear functional on the Banach space $l_\infty$ of bounded  sequences such that for each real sequence $x=\{x_r\}_{r=1}^\infty\in l_\infty$
$$
\liminf_{r\to\infty} x_r\le \Blim_{r\to\infty} x_r\le \limsup_{r\to\infty} x_r
$$
(~we use the customary notation $\Blim\limits_{r\to\infty} x_r$ for the the Banach limit of the sequence $x$~).

In order to justify the notion of $H$-measures, we will need the following result on representation of bilinear functionals.

\begin{lem}\label{Lem5} Let $X,Y$ be locally compact Hausdorff  spaces, and $F(f,g)$ be a bilinear functional on $C_0(X)\times C_0(Y)$ such that for every compact subsets $K_1\subset X$, $K_2\subset Y$
\begin{equation}\label{lc}
|F(f,g)|\le C(K_1,K_2)\|f\|_\infty\|g\|_\infty \ \forall f\in C_0(K_1), g\in C_0(K_2),  \ \mbox{ (continuity), }
\end{equation}
where the constant $C(K_1,K_2)$ depends only on compacts $K_1$, $K_2$, and
\begin{equation}
\label{ln}
F(f,g)\ge 0 \quad \forall f,g\ge 0 \ \mbox{ (nonnegativity)}.
\end{equation}
 Then there exists a unique locally finite nonnegative Radon measure $\mu=\mu(x,y)$ on $X\times Y$ such that
\begin{equation}\label{repr}
F(f,g)=\int_{X\times Y}f(x)g(y)d\mu(x,y).
\end{equation}
\end{lem}
\begin{pf}
First, we consider the case when $X,Y$ are compact sets of Euclidean spaces: $X\subset\R^m$, $Y\subset\R^l$.
In this case the statement of Lemma~\ref{Lem5} was established in \cite[Lemma~1.10]{Ta2}. For completeness we reproduce below the proof. Assuming that $m\ge l$, we may suppose that $X,Y$ are compact subsets of the same Euclidean space: $X,Y\subset\R^m$. We choose a function $K(z)\in C_0(\R^m)$ such that $K(z)\ge 0$, $\supp K\subset B_1\doteq\{ \ z\in\R^m \ | \ |z|\le 1 \ \}$, $\int K(z)dz=1$, and set $K_r(z)=r^mK(rz)$, where $r\in\N$. Obviously, the sequence $K_r(z)$ converges as $r\to\infty$ to the Dirac $\delta$-measure $\delta(z)$ weakly in $\D'(\R^m)$. For $f(x)\in C(\R^m)$ we introduce the averaged functions $f_r(p)=f* K_r(p)=\int_{\R^m} f(x)K_r(p-x)dx$. By the known properties of averaged functions, $f_r\to f$ as $r\to\infty$ uniformly on any compact. This together with the continuity assumption implies that
\begin{equation}\label{l1}
F(f,g)=\lim_{r\to\infty} F_r(f,g),
\end{equation}
where $F_r(f,g)=F(f_r,g_r)$, and the averaged functions
$$
f_r(p)=\int_{\R^m} f(x)K_r(p-x)dx, \quad g_r(q)=\int_{\R^m} g(y)K_r(q-y)dy
$$
are reduced to the sets $X$ and $Y$, respectively. As it follows from the continuity of $F$,
\begin{equation}\label{l2}
F_r(f,g)=F(f_r,g_r)=\int_{\R^m\times\R^m} f(x)g(y)\alpha_r(x,y)dxdy,
\end{equation}
where
$$
\alpha_r(x,y)=F(K_r(p-x),K_r(q-y)).
$$
It is easy to verify that $\alpha_r(x,y)\in C_0(\R^m\times\R^m)$, $\supp \alpha_r\subset X_r\times Y_r$, where
$X_r=X+B_{1/r}$, $Y_r=Y+B_{1/r}$, $r\in\N$, and by $B_\rho$ we denotes the closed ball of radius $\rho$ centered at zero: $B_\rho=\{ \ z\in\R^m \ | \ |z|\le\rho \ \}$. Moreover, by the nonnegativity of $F$ we see that the functionals $F_r$ are also nonnegative: $F_r(f,g)\ge 0$ whenever $f,g\ge 0$, and this readily implies that the kernels $\alpha_r(x,y)\ge 0$. Besides,
$$
\int_{\R^m\times\R^m}\alpha_r(x,y)dxdy=F_r(1,1)\le C,
$$
where $C=C(X,Y)$ is the constant from (\ref{lc}). Therefore, the sequence of nonnegative measures $\mu_r=\alpha_r(x,y)dxdy$
weakly converges as $r\to\infty$ to a finite nonnegative Radon measure $\mu=\mu(x,y)$. Since
$\displaystyle X\times Y=\cap_{r=1}^\infty X_r\times Y_r$, we see that $\supp\mu\subset X\times Y$. For $f\in C(X)$, $g\in C(Y)$ let $\tilde f, \tilde g\in C(\R^m)$ be continuous extensions of these functions on the whole space. Then, in view of (\ref{l1}), (\ref{l2})
\begin{eqnarray*}
F(f,g)=\lim_{r\to\infty}F_r(\tilde f,\tilde g)=\lim_{r\to\infty}\int_{\R^m\times\R^m} \tilde f(x)\tilde g(y)\alpha_r(x,y)dxdy=\\
\int_{X\times Y} \tilde f(x)\tilde g(y) d\mu(x,y)=\int_{X\times Y} f(x)g(y) d\mu(x,y),
\end{eqnarray*}
and representation (\ref{repr}) follows.
Observe that the measure $\mu$ is finite and uniquely defined by (\ref{repr}) because linear combinations of the functions
$f(x)g(y)$ are dense in $C(X\times Y)$. Thus, the proof in the case of compacts $X\subset\R^m$, $Y\subset\R^l$ is complete.

In the case of arbitrary Hausdorff compacts $X$, $Y$, we introduce the set $\mathfrak{A}$, consisting of pairs $(A,B)$ of finite subsets $A\subset C(X)$, $B\subset C(Y)$. The set $\mathfrak{A}$ is ordered by the inclusion order: $\alpha=(A_1,B_1)\le \beta=(A_2,B_2)$ if $A_1\subset A_2$, $B_1\subset B_2$. It is clear, that
for each $\alpha,\beta\in\mathfrak{A}$ there exists $\gamma\in\mathfrak{A}$ such that $\alpha\le\gamma$, $\beta\le\gamma$, that is, $\mathfrak{A}$ is a directed set. Let $\alpha=(A,B)\in\mathfrak{A}$,
$A=\{f_1(x),\ldots,f_m(x)\}\subset C(X)$, $B=\{g_1(y),\ldots,g_l(y)\}\subset C(Y)$, $m,l\in\N$, and let
$F:X\mapsto\R^m$, $G:Y\mapsto\R^l$ be continuous mapping such that $F(x)=(f_1(x),\ldots,f_m(x))$, $G(y)=(g_1(y),\ldots,g_l(y))$, Then $\tilde X=F(X)$, $\tilde Y=G(Y)$ are compact subsets of Euclidean spaces $\R^m$ and $\R^l$, respectively.
We introduce the bilinear functional $F_\alpha(\phi,\psi)$ on $C(\tilde X)\times C(\tilde Y)$, setting
\begin{equation}\label{l3}
F_\alpha(\phi,\psi)=F(\phi(F(x)),\psi(G(y))).
\end{equation}
Clearly, this functional satisfies both the continuity and the nonnegativity conditions. Then, as we have already established, there exists a unique nonnegative Radon measure $\nu_\alpha=\nu_\alpha(p,q)$ on $\tilde X\times\tilde Y$ such that
\begin{equation}\label{l4}
F_\alpha(\phi,\psi)=\int_{\tilde X\times\tilde Y} \phi(p)\psi(q)d\nu_\alpha(p,q).
\end{equation}
Moreover, $\nu_\alpha(\tilde X\times\tilde Y)=F(1,1)\le C$, where $C=C(\tilde X,\tilde Y)$ is the constant from condition (\ref{lc}).
We consider the linear functional
\begin{equation}\label{l5}
\varphi_\alpha(h)=\int \tilde h(p,q)d\nu_\alpha(p,q),
\end{equation}
defined on the subspace $H_\alpha$ of $C(X\times Y)$, consisting of functions $h(x,y)=\tilde h(F(x),G(y))$,
$\tilde h(p,q)\in C(\tilde X\times\tilde Y)$. This functional satisfies the property
\begin{equation}\label{l6}
\varphi_\alpha(h)\le p(h)\doteq C\max_{X\times Y} h^+(x,y), \quad h^+=\max(h,0)
\end{equation}
for all real function $h\in H_\alpha$.
Observe that $p(h)$ is a sub-linear functional on $C(X\times Y)$. Hence, by Hahn-Banach theorem the functional
$\varphi_\alpha$ can be extended to a linear functional $\tilde\varphi_\alpha$ on the whole space $C(X\times Y)$, satisfying estimate (\ref{l6}) for real continuous functions on $X\times Y$. In particular, for each real $h(x,y)\in C(X\times Y)$
$$
-C\max_{X\times Y} h^-(x,y)=-p(-h)\le\tilde\varphi_\alpha(h)\le p(h)=C\max_{X\times Y} h^+(x,y),
$$
which implies, firstly, that $\tilde\varphi_\alpha(h)\ge 0$ whenever $h\ge 0$ and, secondly, that $|\tilde\varphi_\alpha(h)|\le C\max_{X\times Y} |h(x,y)|=C\|h\|_\infty$.
We see that $\tilde\varphi_\alpha$ is a nonnegative continuous functional on $C(X\times Y)$, and $\|\varphi_\alpha\|\le C$. By Riesz-Markov representation theorem there exists a unique nonnegative Radon measure $\mu_\alpha$
on $X\times Y$ such that
\begin{equation}\label{l7}
\tilde\varphi_\alpha(h)=\int_{X\times Y} h(x,y)d\mu_\alpha(x,y),
\end{equation}
and $\mu_\alpha(X\times Y)\le C$.
Observe also that in view of (\ref{l3}), (\ref{l4}), (\ref{l5}), and (\ref{l7})
\begin{equation}\label{l8}
F(f,g)=\varphi_\alpha(f(x)g(y))=\int_{X\times Y} f(x)g(y)d\mu_\alpha(x,y)
\end{equation}
for all $f\in A$, $g\in B$.
Since the space $\M(X\times Y)$ of bounded Radon measures on $X\times Y$ (with the total variation as a norm) is dual to $C(X\times Y)$, then bounded sets in $\M(X\times Y)$ are weakly precompact. Therefore, there exists an accumulation point $\mu$ of a net $\mu_\alpha$, $\alpha\in\mathfrak{A}$
with respect to the weak topology in $\M(X\times Y)$. Let $f(x)\in C(X)$, $g(y)\in C(Y)$, and $\alpha_0=(\{f\},\{g\})\in\mathfrak{A}$. Since $\mu$ is an accumulation point of a net $\mu_\alpha$, then there exists a increasing sequence $\alpha_n=(A_n,B_n)\in\mathfrak{A}$, $n\in\N$, such that $\alpha_n>\alpha_0$ and
in view of (\ref{l8})
$$
F(f,g)=\varphi_{\alpha_n}(f(x)g(y))=\int_{X\times Y} f(x)g(y)d\mu_{\alpha_n}(x,y)\mathop{\to}_{n\to\infty}\int_{X\times Y} f(x)g(y)d\mu(x,y).
$$
This relation implies the desired representation (\ref{repr}) with the finite nonnegative Radon measure $\mu$. Uniqueness of the measure $\mu$ follows again from the density in $C(X\times Y)$ of linear combinations of the functions $f(x)g(y)$.

Now, we consider the general case of locally compact Hausdorff spaces $X,Y$. We introduce the directed set $\mathfrak{K}$ consisting of pairs $\alpha=(K,L)$ of compacts $K\subset X$, $L\subset Y$ and ordered by the inclusion order, i.e., $\alpha=(K,L)\le\alpha_1=(K_1,L_1)$ if $K\subset K_1$, $L\subset L_1$. For each $\alpha=(K,L)\in\mathfrak{K}$ there exist functions $a_\alpha(x)\in C_0(X)$, $b_\alpha(y)\in C_0(Y)$ with the following properties $0\le a_\alpha(x)\le 1$, $0\le b_\alpha(y)\le 1$, and $a_\alpha(x)=b_\alpha(y)=1$ for all $x\in K$, $y\in L$. We denote $X_\alpha=\supp a_\alpha$, $Y_\alpha=\supp b_\alpha$ and define the bilinear functional
$F\alpha:C(X_\alpha)\times C(Y_\alpha)\to\C$ by the identity $F_\alpha(f,g)=F(fa_\alpha,gb_\alpha)$. It is assumed that the functions $(fa_\alpha)(x)$, $(gb_\alpha)(y)$ are extended on the whole spaces $X$, $Y$, being zero outside of $X_\alpha$, $Y_\alpha$, respectively. In particular, these functions have compact supports and the functional $F_\alpha$ is well-defined. Obviously,
$$
F_\alpha(f,g)\le C(X_\alpha,Y_\alpha)\|fa_\alpha\|_\infty\|gb_\alpha\|_\infty\le C_\alpha\|f\|_\infty\|g\|_\infty, \quad C_\alpha=C(X_\alpha,Y_\alpha)
$$
and $F_\alpha(f,g)=F(fa_\alpha,gb_\alpha)\ge 0$ whenever $f,g$ are real and nonnegative. Since $X_\alpha$, $Y_\alpha$ are compact, then, as it was already established above, there exists a unique nonnegative Radon measure $\mu_\alpha$ on $X_\alpha\times Y_\alpha$ such that
$\displaystyle F_\alpha(f,g)=\int_{X_\alpha\times Y_\alpha}f(x)g(y)d\mu_\alpha(x,y).$ This measure may be considered as a Radon measure on the space $X\times Y$ with the support in $X_\alpha\times Y_\alpha$. Then, for every $f=f(x)\in C(X_\alpha)$, $g=g(y)\in C(Y_\alpha)$
\begin{equation}\label{l9}
F_\alpha(f,g)=\int_{X\times Y}f(x)g(y)d\mu_\alpha(x,y).
\end{equation}
Let $K\subset X$, $L\subset Y$ be compact subsets, $\beta=(K,L)\in\mathfrak{K}$. Suppose that $\alpha\in\mathfrak{K}$.
Since the nonnegative function $a_\beta(x)b_\beta(y)\equiv 1$ on $K\times L$, then
\begin{eqnarray}\label{l10}
\mu_\alpha(K\times L)\le \int_{X\times Y} a_\beta(x)b_\beta(y)d\mu_\alpha(y)=F_\alpha(a_\beta,b_\beta)=\nonumber\\
F(a_\alpha a_\beta, b_\alpha b_\beta)\le F(a_\beta, b_\beta)=F_\beta(1,1)\le C_\beta,
\end{eqnarray}
where we use the nonnegativity of $F$, that implies the monotonicity of this functional on the sets of nonnegative functions:
$F(f_1,g_1)\ge F(f_2,g_2)$ for all $f_1,f_2\in C_0(X)$, $g_1,g_2\in C_0(Y)$ such that $0\le f_1(x)\le f_2(x)$, $0\le g_1(y)\le g_2(y)$
(indeed,
$F(f_2,g_2)-F(f_1,g_1)=F(f_2-f_1,g_2)+F(f_1,g_2-g_1)\ge 0$).

In view of estimates (\ref{l10}), the net $\mu_\alpha$, $\alpha\in\mathfrak{K}$ is bounded in locally convex space
$\M_{loc}(X\times Y)$ of locally finite Radon measures (with topology generated by seminorms $p_\alpha(\mu)=|\mu|(K\times L)$, $\alpha=(K,L)\in\mathfrak{K}$, $|\mu|$ stands for the variation of measure $\mu$.

Since the bounded sets of the space $\M_{loc}(X\times Y)$ (which is dual to $C_0(X\times Y)$) are compact, there exists a weak accumulation point $\mu\in \M_{loc}(X\times Y)$ of the net $\mu_\alpha$, $\alpha\in\mathfrak{K}$. Since $\mu_\alpha\ge 0$ for all $\alpha\in\mathfrak{K}$, we claim that $\mu\ge 0$. Let $f(x)\in C_0(X)$, $g(y)\in C_0(Y)$, and $\alpha_0=(\supp f,\supp g)\in\mathfrak{K}$. Since $\mu$ is an accumulation point of the net $\mu_\alpha$, $\alpha\in\mathfrak{K}$, there exists a increasing sequence $\alpha_n=(K_n,L_n)\in\mathfrak{K}$, $n\in\N$, such that $\alpha_n>\alpha_0$ and
$$
F(f,g)=\int_{X\times Y} f(x)g(y)d\mu_{\alpha_n}(x,y)\mathop{\to}_{n\to\infty}\int_{X\times Y} f(x)g(y)d\mu(x,y).
$$
This implies representation (\ref{repr}) and conclude the proof.
\end{pf}

The following statement, analogous to the assertions of Propositions~\ref{Pro1},\ref{Pro2}, holds.

\begin{prop}\label{Pro5} There exists a family of Radon measures $\mu=\{\mu^{\alpha\beta}\}_{\alpha,\beta=1}^N$ on $\Omega\times\S$
such that for all  $\Phi_1(x),\Phi_2(x)\in C_0(\Omega)$, $\psi(\xi)\in\A$, and $\alpha,\beta=1,\ldots,N$
\begin{equation}\label{Hmnew}
\<\mu^{\alpha\beta}(x,\eta),\Phi_1(x)\overline{\Phi_2(x)}\hat\psi(\eta)\>= \Blim_{r\to\infty}\int_{\R^n}
F(U_r^\alpha\Phi_1)(\xi)\overline{F(U_r^\beta\Phi_2)(\xi)} \psi(\xi)d\xi.
\end{equation}
The matrix-valued measure $\mu$ is Hermitian and positive semi-definite, i.e., for every  $\zeta=(\zeta_1,\ldots,\zeta_N)\in\C^n$
$$
\mu\zeta\cdot\zeta=
\sum_{\alpha,\beta=1}^N\mu^{\alpha\beta}\zeta_\alpha\overline{\zeta_\beta}\ge 0.
$$
\end{prop}
\begin{pf}
Denote for $\Phi_1(x),\Phi_2(x)\in C_0(\Omega)$, $\psi(\xi)\in A$
\begin{equation}\label{1}
I^{\alpha\beta}(\Phi_1,\Phi_2,\psi)=\Blim_{r\to\infty}\int_{\R^n}
F(\Phi_1U_r^\alpha)(\xi)\overline{F(\Phi_2U_r^\beta)(\xi)}\psi(\xi)d\xi
\end{equation}
and observe that, by the Buniakovskii inequality and the Plancherel identity,
\begin{equation}\label{2}
|I^{\alpha\beta}|\le
\|\Phi_1\|_\infty\|\Phi_2\|_\infty\|\psi\|_\infty\cdot\limsup_{r\to\infty}
\left[\|U_r^\alpha\|_{L^2(K)}\|U_r^\beta\|_{L^2(K)}\right],
\end{equation}
where $K\subset\Omega$ is a compact containing supports of $\Phi_1$ and $\Phi_2$.
In view of the weak convergence of sequences $U_r^\alpha$ in $L^2(K)$ these sequences are bounded in $L^2(K)$.
Therefore, for some constant $C_K$ we have $\|U_r^\alpha\|_{L^2(K)}^2\le C_K$ for all $r\in\N$,
$\alpha=1,\ldots,N$. Then, it follows from (\ref{2}) that
\begin{equation}
\label{3}
|I^{\alpha\beta}(\Phi_1,\Phi_2,\psi)|\le C_K\|\Phi_1\|_\infty\|\Phi_2\|_\infty\|\psi\|_\infty
\end{equation}
with $K=\supp\Phi_1\cup\supp\Phi_2$.
If $\psi(\xi)\in A_0$, then by Lemma~\ref{Lem1} the operator $B_\psi A_{\Phi_1}$ is compact in $L^2(\R^n)$. Hence, the sequences $B_\psi A_{\Phi_1}(U_r^\alpha)=B_\psi A_{\Phi_1}(U_r^\alpha\chi_K)\to 0$ in $L^2(\R^n)$. Here $\chi_K(x)$
is the indicator function of the compact $K=\supp\Phi_1$. We see that for all $\alpha=1,\ldots,N$
$$
F(\Phi_1U_r^\alpha)(\xi)\psi(\xi)=F(B_\psi A_{\Phi_1}(U_r^\alpha))(\xi)\mathop{\to}_{r\to\infty} 0 \ \mbox{ in } L^2(\R^n),
$$
and, therefore, for all $\alpha,\beta=1,\ldots,N$
$$
\lim_{r\to\infty}\int_{\R^n}
F(\Phi_1U_r^\alpha)(\xi)\overline{F(\Phi_2U_r^\beta)(\xi)}\psi(\xi)d\xi=0.
$$
In view of (\ref{1})  $I^{\alpha\beta}(\Phi_1,\Phi_2,\psi)=0$ for $\Phi_1(x),\Phi_2(x)\in C_0(\Omega)$ and all $\psi(\xi)\in A_0$. We see that the linear with respect to $\psi$ functional $I^{\alpha\beta}(\Phi_1,\Phi_2,\psi)$
is well-defined on factor-algebra $\A=A/A_0$ and, in view of (\ref{3}), for all $\psi_0\in A_0$
$$
|I^{\alpha\beta}(\Phi_1,\Phi_2,\psi)|=|I^{\alpha\beta}(\Phi_1,\Phi_2,\psi-\psi_0)|\le C_K\|\Phi_1\|_\infty\|\Phi_2\|_\infty\|\psi-\psi_0\|_\infty.
$$
Therefore,
\begin{equation}\label{4}
I^{\alpha\beta}(\Phi_1,\Phi_2,\psi)|\le C_K\|\Phi_1\|_\infty\|\Phi_2\|_\infty \inf_{\psi_0\in A_0}\|\psi-\psi_0\|_\infty=C_K\|\Phi_1\|_\infty\|\Phi_2\|_\infty\|\psi\|_{\A},
\end{equation}
where
$$
\|\psi\|_{\A}=\inf_{\psi_0\in A_0}\|\psi-\psi_0\|_\infty=\esslimsup_{|\xi|\to\infty}|\psi(\xi)|
$$
is the factor-norm of $[\psi]$ in $\A$.
Now, we observe that
\begin{equation}
\label{5}
\int_{\R^n}
F(\Phi_1U_r^\alpha)(\xi)\overline{F(\Phi_2U_r^\beta)(\xi)}\psi(\xi)d\xi=(A_\psi(\Phi_1U_r^\alpha),\Phi_2U_r^\beta)_2,
\end{equation}
where $(\cdot,\cdot)_2$ is the scalar product in $L^2=L^2(\R^n)$. Let $\omega(x)\in C_0(\R^n)$ be a
function such that $\omega(x)\equiv 1$ on $\supp\Phi_1$. Then
\begin{equation}\label{6}
A_\psi(\Phi_1U_r^\alpha)=A_\psi B_{\Phi_1}(\omega U_r^\alpha)=B_{\Phi_1}A_\psi(\omega U_r^\alpha)+[A_\psi,B_{\Phi_1}](\omega U_r^\alpha).
\end{equation}
By the definition of algebra $\A$, the operator $[A_\psi,B_{\Phi_1}]$
is compact on $L^2$ and
since $\omega U_r^\alpha\rightharpoonup 0$ as $r\to\infty$ weakly in $L^2$, we claim that $[A_\psi,B_{\Phi_1}](\omega U_r^\alpha)\to 0$ as $r\to\infty$ strongly
in $L^2$. Since the sequence $\Phi_2U_r^\beta$ is bounded in $L^2$, we conclude that
$([A_\psi,B_{\Phi_1}](\omega U_r^\alpha),\Phi_2 U_r^\beta))_2\to 0$ as $r\to\infty$. It follows from this limit
relation and (\ref{5}), (\ref{6}) that
$$
I^{\alpha\beta}(\Phi_1,\Phi_2,\psi)=\Blim_{r\to\infty} (B_{\Phi_1}A_\psi(\omega U_r^\alpha),\Phi_2 U_r^\beta))_2=
\Blim_{r\to\infty}\int_{\R^n} \Phi_1(x)\overline{\Phi_2(x)}A_\psi(\omega U_r^\alpha)(x)\overline{U_r^\beta(x)}dx.
$$
We claim that
$$
I^{\alpha\beta}(\Phi_1,\Phi_2,\psi)=\tilde I^{\alpha\beta}(\Phi_1\overline{\Phi_2},\hat\psi),
$$
where $$\tilde I^{\alpha\beta}(\Phi,\hat\psi)=\Blim_{r\to\infty}\int_{\R^n} \Phi(x)A_\psi(\omega U_r^\alpha)(x)\overline{U_r^\beta(x)}dx$$ is a bilinear functional on $C_0(\Omega)\times C(\S)$ for
each $\alpha,\beta=1,\ldots,N$ ($\omega(x)\in C_0(\R^n)$ is an arbitrary function equalled $1$ on support of $\Phi(x)$), and $\hat\psi(\eta)$ being the Gelfand transform of $\psi(\xi)$. Taking in the above relation $\Phi_1=\Phi(x)/\sqrt{|\Phi(x)|}$ (we set
$\Phi_1(x)=0$ if $\Phi(x)=0$), $\Phi_2=\sqrt{|\Phi(x)|}$, where $\Phi(x)\in C_0(\Omega)$, we find with
the help of (\ref{4}) that
\begin{eqnarray*}
|\tilde I^{\alpha\beta}(\Phi,\hat\psi)|=|I^{\alpha\beta}(\Phi_1,\Phi_2,\psi)|\le
C_K\|\Phi_1\|_\infty\|\Phi_2\|_\infty\|\psi\|_{\A} \\ = C_K\|\Phi\|_\infty\|\psi\|_\A=C_K\|\Phi\|_\infty\|\hat\psi\|_\infty, \
K=\supp\Phi.
\end{eqnarray*}
This estimate shows that the functionals $\tilde I^{\alpha\beta}(\Phi,\hat\psi)$ are continuous on
$C_0(\Omega)\times C(\S)$. Now, we observe that for nonnegative $\Phi(x)$ and $\hat\psi(\eta)$ the matrix
$\tilde I\doteq\{\tilde I^{\alpha\beta}(\Phi,\tilde\psi)\}_{\alpha,\beta=1}^N$ is Hermitian and positive
definite. First, we remark that by Lemma~\ref{Lem2} \ $\hat\psi(\eta)\ge 0$ if and only if $\psi(\xi)\ge 0$. Taking $\Phi_1(x)=\Phi_2(x)=\sqrt{\Phi(x)}$, we find
\begin{equation}\label{7}
\tilde I^{\alpha\beta}(\Phi,\hat\psi)=I^{\alpha\beta}(\Phi_1,\Phi_1,\psi)=
\Blim_{r\to\infty}\int_{\R^n} F(\Phi_1U_r^\alpha)(\xi)\overline{F(\Phi_1U_r^\beta)(\xi)}
\psi(\xi)d\xi.
\end{equation}
For $\zeta=(\zeta_1,\ldots,\zeta_N)\in\C^N$ we have, in view of (\ref{7}),
$$
\tilde I\zeta\cdot\zeta=\sum_{\alpha,\beta=1}^N\tilde
I^{\alpha\beta}(\Phi,\hat\psi)\zeta_\alpha\overline{\zeta_\beta}=\Blim_{r\to\infty}\int_{\R^n}
|F(\Phi_1V_r)(\xi)|^2\psi(\xi)d\xi\ge 0,
$$
where $V_r(x)=\sum\limits_{\alpha=1}^N U_r^\alpha\zeta_\alpha$. The above relation proves that the
matrix $\tilde I$ is Hermitian and positive definite.

We see that for any $\zeta\in\C^n$ the bilinear functional $\tilde I(\Phi,\hat\psi)\zeta\cdot\zeta$ is
continuous on $C_0(\Omega)\times C(\S)$ and nonnegative, that is, $\tilde
I(\Phi,\hat\psi)\zeta\cdot\zeta\ge 0$ whenever $\Phi(x)\ge 0$, $\hat\psi(\eta)\ge 0$. By Lemma~\ref{Lem5} such a functional is represented by integration over
some unique locally finite non-negative Radon measure $\mu=\mu_\zeta(x,\eta)\in M_{loc}(\Omega\times\S)$:
$$
\tilde I(\Phi,\hat\psi)\zeta\cdot\zeta=\int_{\Omega\times\S}\Phi(x)\hat\psi(\eta)d\mu_\zeta(x,\eta).
$$
As a function of the vector $\zeta$, $\mu_\zeta$ is a measure valued Hermitian form. Therefore,
\begin{equation}\label{8}
\mu_\zeta=\sum_{\alpha,\beta=1}^N \mu^{\alpha\beta}\zeta_\alpha\overline{\zeta_\beta}
\end{equation}
with measure valued coefficients $\mu^{\alpha\beta}\in M_{loc}(\Omega\times\S)$, which can be
expressed as follows
$$
\mu^{\alpha\beta}=[\mu_{e_\alpha+e_\beta}+i\mu_{e_\alpha+ie_\beta}]/2-(1+i)(\mu_{e_\alpha}+\mu_{e_\beta})/2,
$$
where $e_1,\ldots,e_N$ is the standard basis in $\C^N$, and $i^2=-1$.

By (\ref{8})
$$
\tilde
I(\Phi,\hat\psi)\zeta\cdot\zeta=\sum_{\alpha,\beta=1}^N\<\mu^{\alpha\beta},\Phi(x)\hat\psi(\eta)\>\zeta_\alpha\overline{\zeta_\beta}
$$
and since
$$\tilde I(\Phi,\hat\psi)\zeta\cdot\zeta=\sum_{\alpha,\beta=1}^N \tilde I^{\alpha\beta}(\Phi,\hat\psi)\zeta_\alpha\overline{\zeta_\beta},$$
then, comparing the coefficients, we find that
\begin{equation}\label{9}
\<\mu^{\alpha\beta},\Phi(x)\hat\psi(\eta)\>=\tilde I^{\alpha\beta}(\Phi,\hat\psi).
\end{equation}
In particular,
$$
\<\mu^{\alpha\beta},\Phi_1(x)\overline{\Phi_2(x)}\hat\psi(\eta)\>=I^{\alpha\beta}(\Phi_1,\Phi_2,\psi)=
\Blim_{r\to\infty}\int_{\R^n}
F(\Phi_1U_r^\alpha)(\xi)\overline{F(\Phi_2U_r^\beta)(\xi)} \psi(\xi)d\xi.
$$
To complete the proof, observe that for each $\zeta\in\C^N$ the
measure
$$
\sum_{\alpha,\beta=1}^N\mu^{\alpha\beta}\zeta_\alpha\overline{\zeta_\beta}=\mu_\zeta\ge 0.
$$
Hence, $\mu$ is Hermitian and positive definite.
\end{pf}

The usage of generalized Banach limit instead of extraction of a subsequence is connected with the fact that the algebra $\A$ is not separable. Therefore, the extraction of a subsequence of $U_r$ such that relation (\ref{Hmnew}) holds, with replacement of the Banach limit to the usual one, is not always possible.
Certainly, the $H$-measure $\mu$ depend on the choice of the generalized Banach limit (this resembles the dependence of the Tartar $H$-measure on the choice of a subsequence).

If $X$ is a subspace of $\R^n$ and $p_X:\S\to S_X$ is the projection defined before Proposition~\ref{Pro3} above, then the image
of the measures $\mu^{\alpha\beta}$ under the map $(x,\eta)\to (x,p_X(\eta))$ is exactly the ultraparabolic
$H$-measure corresponding to the subspace $X$.

Evidently, if the sequence $U_r$ converges as $r\to\infty$ to the zero vector strongly in $L^2_{loc}(\Omega,\C^N)$, then $H$-measure is trivial: $\mu=0$. Conversely, if $\mu=0$ then for each
$\Phi(x)\in C_0(\Omega)$
$$
\Blim_{r\to\infty}\int_\Omega |U_r(x)\Phi(x)|^2dx=\sum_{\alpha=1}^N\Blim_{r\to\infty} \int_\Omega |F(U_r^\alpha\Phi)(\xi)|^2d\xi=\sum_{\alpha=1}^N\<\mu^{\alpha\alpha}(x,\xi),|\Phi(x)|^2\>=0.
$$
This implies that
\begin{equation}\label{32}
\liminf_{r\to\infty}\int_\Omega |U_r(x)\Phi(x)|^2dx=0 \ \forall\Phi(x)\in C_0(\Omega).
\end{equation}
We can choose the sequence of real nonnegative functions $\Phi_k(x)\in C_0(\Omega)$ such that
$\Phi_{k+1}(x)\ge\Phi_k(x)$ for all $k\in\N$, and $\lim\limits_{k\to\infty}\Phi_k(x)=1$ for all $x\in\Omega$. It follows from (\ref{32}) that there exists a strictly increasing sequence $r_k\in\N$ such that $\displaystyle\int_\Omega |U_{r_k}(x)\Phi_k(x)|^2dx<1/k$. Then the subsequence
$$
U_{r_k}(x)\mathop{\to}_{k\to\infty} 0 \ \mbox{ in } L^2_{loc}(\Omega,\C^N).
$$

\medskip
Let $\mu=\{\mu^{\alpha\beta}\}_{\alpha,\beta=1}^N$ be an $H$-measure corresponding to a sequence $U_r=\{U_r^\alpha\}_{\alpha=1}^N\in L^2_{loc}(\Omega,\C^N)$
We define $\mu_0={\rm
Tr}\mu=\sum_{\alpha=1}^N \mu^{\alpha\alpha}$. As follows from
Proposition~\ref{Pro5}, \ $\mu_0$ is a locally finite non-negative Radon measure on
$\Omega\times\S$. We assume that this measure is extended on
$\sigma$-algebra of $\mu_0$-measurable sets, and in particular that
this measure is complete.

\begin{lem}\label{Lem6}
The $H$-measure $\mu$ is absolutely continuous with respect to $\mu_0$, more precisely,
$\mu=H(x,\eta)\mu_0$, where $H(x,\eta)=\{h^{\alpha\beta}(x,\eta)\}_{\alpha,\beta=1}^N$ is a bounded
$\mu_0$-measurable function taking values in the cone of nonnegative definite Hermitian $N\times N$
matrices, moreover $|h^{\alpha\beta}(x,\eta)|\le 1$.
\end{lem}
\begin{pf}
Remark firstly that $\mu^{\alpha\alpha}\le\mu_0$ for all $\alpha=1,\ldots,N$. Now,
suppose that $\alpha,\beta\in\{1,\ldots,N\}$, $\alpha\not=\beta$. By Proposition~\ref{Pro5} for any
compact set $B\subset\Omega\times\S$ the matrix
$$
\left(\begin{array}{cc} \mu^{\alpha\alpha}(B) &
\mu^{\alpha\beta}(B) \\[5pt]
\overline{\mu^{\alpha\beta}(B)} & \mu^{\beta\beta}(B)
\end{array}\right)
$$
is nonnegative definite; in particular,
$$
|\mu^{\alpha\beta}(B)|\le\left(\mu^{\alpha\alpha}(B)\mu^{\beta\beta}(B)\right)^{1/2}\le\mu_0(B).
$$
By regularity of measures $\mu^{\alpha\beta}$ and $\mu_0$ this estimate is satisfied for all Borel sets
$B$. This easily implies the inequality $\Var\mu^{\alpha\beta}\le\mu_0$. In particular, the measures
$\mu^{\alpha\beta}$ are absolutely continuous with respect to $\mu_0$, and by the Radon-Nykodim theorem
$\mu^{\alpha\beta}=h^{\alpha\beta}(x,\eta)\mu_0$, where the densities $h^{\alpha\beta}(x,\eta)$ are
$\mu_0$-measurable and, as follows from the inequalities $\Var\mu^{\alpha\beta}\le\mu_0$,
$|h^{\alpha\beta}(x,\eta)|\le 1$ $\mu_0$-a.e. on $\Omega\times\S$. We denote by $H(x,\eta)$ the matrix
with components $h^{\alpha\beta}(x,\eta)$. Recall that the $H$-measure $\mu$ is nonnegative definite. This
means that for all $\zeta\in\C^N$
\begin{equation}\label{CC9}
\mu\zeta\cdot\zeta=(H(x,\eta)\zeta\cdot\zeta)\mu_0\ge 0.
\end{equation}
Hence $H(x,\eta)\zeta\cdot\zeta\ge 0$ for $\mu_0$-a.e. $(x,\eta)\in\Omega\times\S$. Choose a countable
dense set $E\subset \C^N$. Since $E$ is countable, then it follows from (\ref{CC9}) that for a set
$(x,\eta)\in\Omega\times\S$ of full $\mu_0$-measure $H(x,\eta)\zeta\cdot\zeta\ge 0$ $\forall\zeta\in
E$, and since $E$ is dense we conclude that actually $H(x,\eta)\zeta\cdot\zeta\ge 0$ for all
$\zeta\in\C^N$. Thus, the matrix $H(x,\eta)$ is Hermitian and nonnegative definite for $\mu_0$-a.e.
$(x,\eta)$. After an appropriate correction on a set of null $\mu_0$-measure, we can assume that the
above property is satisfied for all $(x,\eta)\in\Omega\times\S$, and also
$|h^{\alpha\beta}(x,\eta)|\le 1$ for all $(x,\eta)\in\Omega\times\S$, $\alpha,\beta=1,\ldots,N$. The
proof is complete.
\end{pf}

Now we assume that for all $\Phi(x)\in C_0^\infty(\Omega)$ the sequence
\begin{equation}\label{33}
\sum_{\alpha=1}^N\sum_{k=1}^M c_{\alpha k}(x)p_{\alpha k}(\partial/\partial x)(\Phi(x)U^\alpha_r(x))\mathop{\to}_{r\to\infty} 0 \ \mbox{ in } L^2_{loc}(\Omega),
\end{equation}
where $p_{\alpha k}(\partial/\partial x)$ denotes the pseudo-differential operator with symbol ${p_{\alpha k}(\xi)\in A}$ and $c_{\alpha k}(x)\in C(\Omega)$. Then the $H$-measure corresponding to the sequence $U_r$ satisfy the following localization property.

\begin{thm}\label{Th2} For each $\beta=1,\ldots,N$
$$
\sum_{\alpha=1}^N\sum_{k=1}^M c_{\alpha k}(x)\widehat{p_{\alpha k}}(\eta)\mu^{\alpha\beta}(x,\eta)=0.
$$
\end{thm}
\begin{pf}
In view of (\ref{33}) for all $\Phi_1(x),\Phi(x)\in C_0^\infty(\Omega)$ a sequence
\begin{equation}\label{34}
\sum_{\alpha=1}^N\sum_{k=1}^M \Phi(x)c_{\alpha k}(x)p_{\alpha k}(\partial/\partial x)(\Phi_1(x)U^\alpha_r(x))\mathop{\to}_{r\to\infty} 0 \ \mbox{ in } L^2(\R^n).
\end{equation}
Since $p_{\alpha k}\in A$, then the commutator $[A_{p_{\alpha k}},B_{\Phi c_{\alpha k}}]$ is a compact operator in $L^2(\R^n)$. Therefore, this operators transform the weakly convergent sequence $\Phi_1(x)U^\alpha_r(x)$ to the strongly convergent ones, which implies that
\begin{equation}\label{35}
\sum_{\alpha=1}^N\sum_{k=1}^M [A_{p_{\alpha k}},B_{\Phi c_{\alpha k}}](\Phi_1(x)U^\alpha_r(x))\mathop{\to}_{r\to\infty} 0 \ \mbox{ in } L^2(\R^n).
\end{equation}
Putting relations (\ref{34}), (\ref{35}) together, we find
$$
\sum_{\alpha=1}^N\sum_{k=1}^M p_{\alpha k}(\partial/\partial x)(\Phi(x)\Phi_1(x)c_{\alpha k}(x)U^\alpha_r(x))
\mathop{\to}_{r\to\infty} 0 \ \mbox{ in } L^2(\R^n).
$$
Taking $\Phi(x)$ in such a way that $\Phi(x)=1$ on $\supp\Phi_1$, we arrive at the relation
$$
\sum_{\alpha=1}^N\sum_{k=1}^M p_{\alpha k}(\partial/\partial x)(\Phi_1(x)c_{\alpha k}(x)U^\alpha_r(x))
\mathop{\to}_{r\to\infty} 0 \ \mbox{ in } L^2(\R^n).
$$
Applying the Fourier transformation, we obtain
\begin{equation}\label{36}
\sum_{\alpha=1}^N\sum_{k=1}^M p_{\alpha k}(\xi)F(\Phi_1c_{\alpha k}U^\alpha_r)(\xi)\mathop{\to}_{r\to\infty} 0 \ \mbox{ in } L^2(\R^n).
\end{equation}
We multiply (\ref{36}) by the bounded sequence $\overline{F(\Phi_2U_r^\beta)(\xi)}\psi(\xi)$, where
$\Phi_2(x)\in C_0^\infty(\Omega)$, $\psi(\xi)\in A$, and $1\le\beta\le N$. Integrating over $\xi\in\R^n$, we
arrive at the relation
$$
\lim_{r\to\infty}\sum_{\alpha=1}^N\sum_{k=1}^M\int_{\R^n} p_{\alpha k}(\xi)F(\Phi_1c_{\alpha k}U^\alpha_r)(\xi)\overline{F(\Phi_2U_r^\beta)(\xi)}\psi(\xi)d\xi=0.
$$
On the other hand, by the definition of $H$-measure this limit coincides with
\begin{eqnarray*}
\Blim_{r\to\infty}\sum_{\alpha=1}^N\sum_{k=1}^M\int_{\R^n} p_{\alpha k}(\xi)F(\Phi_1c_{\alpha k}U^\alpha_r)(\xi)\overline{F(\Phi_2U_r^\beta)(\xi)}\psi(\xi)d\xi= \\
\sum_{\alpha=1}^N\sum_{k=1}^M\Blim_{r\to\infty}\int_{\R^n} p_{\alpha k}(\xi)F(\Phi_1c_{\alpha k}U^\alpha_r)(\xi)\overline{F(\Phi_2U_r^\beta)(\xi)}\psi(\xi)d\xi= \\
\sum_{\alpha=1}^N\sum_{k=1}^M\<\mu^{\alpha\beta},c_{\alpha k}(x)\Phi_1(x)\overline{\Phi_2(x)}\widehat{p_{\alpha k}}(\eta)\widehat{\psi}(\eta)\>.
\end{eqnarray*}
Hence,
$$
\sum_{\alpha=1}^N\sum_{k=1}^M\<\mu^{\alpha\beta},c_{\alpha k}(x)\Phi_1(x)\overline{\Phi_2(x)}\widehat{p_{\alpha k}}(\eta)\widehat{\psi}(\eta)\>=0.
$$
This relation can be written in the form
\begin{equation}\label{37}
\left\<\sum_{\alpha=1}^N\sum_{k=1}^M c_{\alpha k}(x)\widehat{p_{\alpha k}}(\eta)\mu^{\alpha\beta}(x,\eta),
\Phi_1(x)\overline{\Phi_2(x)}\widehat{\psi}(\eta)\right\>=0.
\end{equation}
Since the test functions $\Phi_1(x),\Phi_2(x)\in C_0^\infty(\Omega)$, $\widehat{\psi}(\eta)\in C(\S)$ are arbitrary, the statement of Theorem~\ref{Th2} follows from (\ref{37}).
\end{pf}

\section{Compensated compactness}

Assume that a sequence $u_r(x)\in L^2_{loc}(\Omega,\C^N)$ weakly converges to a vector-function $u(x)$ as $r\to\infty$ while for each $\Phi(x)\in C_0^\infty(\Omega)$ the sequences
\begin{equation}\label{38}
\sum_{\alpha=1}^N\sum_{k=1}^M c_{s\alpha k}(x)p_{s\alpha k}(\partial/\partial x)(\Phi(x)u^\alpha_r(x)) \ \mbox{ are precompact in } L^2_{loc}(\Omega),
\end{equation}
where $p_{s\alpha k}(\partial/\partial x)$ are pseudo-differential operators with symbols ${p_{\alpha k}(\xi)\in A}$, $c_{s\alpha k}(x)\in C(\Omega)$, and $s=1,\ldots,m$. Introduce the set
$$
\Lambda=\Lambda(x)=\Bigl\{ \ \lambda\in\C^N \ | \ \exists\eta\in\S :  \\ \sum_{\alpha=1}^N\sum_{k=1}^M c_{s\alpha k}(x)\widehat{p_{s\alpha k}}(\eta)\lambda_\alpha=0 \ \forall s=1,\ldots,m \ \Bigr\}.
$$
Now, suppose that
$$
q(x,u)=Q(x)u\cdot u=\sum_{\alpha,\beta=1}^N q_{\alpha\beta}(x)u_\alpha \overline{u_\beta}
$$
is an Hermitian form with the matrix $Q(x)$ of coefficients $q_{\alpha\beta}(x)\in C(\Omega)$.

Let the sequence $q(x,u_r)\rightharpoonup v$ as $r\to\infty$ weakly in
$M_{loc}(\Omega)$. The following theorem is analogous to Theorem~\ref{Th1}.

\medskip
\begin{thm}\label{Th3} If $q(x,\lambda)\ge 0$ for all $\lambda\in\Lambda(x)$,
$x\in\Omega$, then ${q(x,u(x))\le v}$ ( in the sense of measures ).
\end{thm}
\begin{pf}
Let $\mu=\{\mu^{\alpha\beta}\}_{\alpha,\beta=1}^N$ be the $H$-measure corresponding to the sequence $U_r=u_r-u$. By Lemma~\ref{Lem6}
this $H$-measure admits the representation $\mu=H(x,\eta)\mu_0$, where $\mu_0=\mathrm{Tr}\,\mu\ge 0$ and
$H(x,\eta)$ is a $\mu_0$-measurable map from $\Omega\times\S$ into the cone of nonnegative definite $N\times N$ Hermitian matrices. As readily follows from (\ref{38}), for each $\Phi(x)\in C_0^\infty(\Omega)$, $s=1,\ldots,m$
$$
\sum_{\alpha=1}^N\sum_{k=1}^M c_{s\alpha k}(x)p_{s\alpha k}(\partial/\partial x)(\Phi(x)U^\alpha_r(x))\mathop{\to}_{r\to\infty} 0  \ \mbox{ strongly in } L^2_{loc}(\Omega).
$$
By Theorem~\ref{Th2} for all $s=1,\ldots,m$ and $\beta=1,\ldots,N$
$$
\sum_{\alpha=1}^N\sum_{k=1}^M c_{s\alpha k}(x)\widehat{p_{s\alpha k}}(\eta)H^{\alpha\beta}(x,\eta)\mu_0=
\sum_{\alpha=1}^N\sum_{k=1}^M c_{s\alpha k}(x)\widehat{p_{s\alpha k}}(\eta)\mu^{\alpha\beta}(x,\eta)=0.
$$
This implies that for $\mu_0$-a.e. $(x,\eta)$ the image $\Im H(x,\eta)\subset\Lambda(x)$. Since the matrix $H(x,\eta)\ge 0$, there exists a unique Hermitian matrix $R=R(x,\eta)=(H(x,\eta))^{1/2}$ such that $R\ge 0$ and $H=R^2$.
By the known properties of Hermitian matrices $\ker R = \ker H$, which readily implies that also
$\Im R = \Im H$. In particular, $\Im R(x,\eta)\subset\Lambda(x)$ for $\mu_0$-a.e. $(x,\eta)\in\Omega\times\S$.
Let $\Phi(x)\in C_0(\Omega)$ be a real test function. Then
\begin{eqnarray}\label{39}
\Blim_{r\to\infty}\int (\Phi(x))^2q(x,U_r(x))=\sum_{\alpha,\beta=1}^N\Blim_{r\to\infty}
\int q_{\alpha\beta}(x)\Phi(x)U_r^\alpha(x)\overline{\Phi(x)U_r^\beta(x)}dx= \nonumber\\
\sum_{\alpha,\beta=1}^N\Blim_{r\to\infty}
\int F(q_{\alpha\beta}\Phi U_r^\alpha)(\xi)\overline{F(\Phi U_r^\beta)(\xi)}d\xi=
\sum_{\alpha,\beta=1}^N \<\mu^{\alpha\beta},q_{\alpha\beta}(x)(\Phi(x))^2\>=\nonumber\\
\int_{\Omega\times\S}(\Phi(x))^2\sum_{\alpha,\beta=1}^Nq_{\alpha\beta}(x)h^{\alpha\beta}(x,\eta)d\mu_0(x,\eta).
\end{eqnarray}
Since $H=R^2$ then
$$h^{\alpha\beta}=\sum_{k=1}^N r_{\alpha k}\overline{r_{\beta k}} \quad \forall \alpha,\beta=1,\ldots,N,$$ where
$r_{ij}=r_{ij}(x,\eta)$, $i,j=1,\ldots,N$ are components of the matrix $R$. Therefore,
\begin{equation}\label{40}
\sum_{\alpha,\beta=1}^N q_{\alpha\beta}(x)h^{\alpha\beta}(x,\eta)=\sum_{k=1}^N \sum_{\alpha,\beta=1}^N q_{\alpha\beta}(x)r_{\alpha k}\overline{r_{\beta k}}=\sum_{k=1}^N q(x,Re_k),
\end{equation}
where $e_k$, $k=1,\ldots,N$, is the standard basis in $\C^N$. Since $R(x,\eta)e_k\in\Im R(x,\eta)\subset\Lambda(x)$ for $\mu_0$-a.e. $(x,\eta)\in\Omega\times\S$, then $q(x,R(x,\eta)e_k)\ge 0$  for $\mu_0$-a.e. $(x,\eta)$ and it follows
from (\ref{39}), (\ref{40}) that
\begin{equation}\label{41}
\Blim_{r\to\infty}\int (\Phi(x))^2q(x,U_r(x))\ge 0.
\end{equation}
for all real $\Phi(x)\in C_0(\Omega)$. In view of the weak convergence $u_r\rightharpoonup u$, $q(x,u_r)\rightharpoonup v$ as $r\to\infty$,
$$
q(x,U_r(x))=q(x,u_r(x))+q(x,u(x))-2\mathop{\rm Re } (Q(x)u_r\cdot u)\rightharpoonup v-q(x,u(x)).
$$
in $\M_{loc}(\Omega)$.
Now, it follows from (\ref{41}) that
$$
\<v-q(x,u(x))dx,(\Phi(x))^2\>=\lim_{r\to\infty}\int (\Phi(x))^2q(x,U_r(x))\ge 0
$$
and since the real test function $\Phi(x)$ is arbitrary, $v\ge q(x,u(x))$. The proof is complete.
\end{pf}

\begin{cor}\label{Cor2} Suppose that $q(x,\lambda) = 0$ for all $\lambda\in\Lambda(x)$. Then $v = q(x,u(x))$, that is, the functional $u\to q(x,u)$ is weakly continuous.
\end{cor}
\begin{pf}
Applying Theorem~\ref{Th3} to the quadratic forms $\pm q(x,u)$, we obtain the inequalities $\pm v \ge\pm q(x,u(x))$, which readily imply that $v = q(x,u(x))$.
\end{pf}

\subsection{The case of second order differential constraints}

Now we assume that the sequences
\begin{equation}\label{42}
\sum_{\alpha=1}^N\sum_{k=1}^n \partial_{x_k}(a_{s\alpha k}u_{\alpha r})
+\sum_{\alpha=1}^N\sum_{k,l=1}^n \partial_{x_kx_l}(b_{s\alpha kl}u_{\alpha r}), \quad s=1,\ldots,m,
\end{equation}
are pre-compact in the Sobolev space $H^{-1}_{loc}(\Omega)\doteq W^{-1}_{2,loc}(\Omega)$, where the coefficients $a_{s\alpha k}=a_{s\alpha k}(x)$, $b_{s\alpha kl}=b_{s\alpha kl}(x)$ are supposed to be (generally -- complex-valued) continuous functions on $\Omega$, and
$b_{s\alpha lk}=b_{s\alpha kl}$, $s=1,\ldots,m$, $\alpha=1,\ldots,N$, $k,l=1,\ldots,n$.

We denote by $A_{s\alpha}=A_{s\alpha}(x)$ the vector $\{a_{s\alpha k}\}_{k=1}^n\in\C^n$ and by $B_{s\alpha}=B_{s\alpha}(x)$ the symmetric matrices with components $\{b_{s\alpha kl}\}_{k,l=1}^n$. Let $X_s$
be the maximal linear subspace of $\R^n$ contained in $\R^n\cap\ker B_{s\alpha}(x)$ for all $\alpha=1,\ldots,N$ and $x\in\Omega$.
The following statement easily follows from the definition of the subspace $X_s$:

\begin{lem}\label{Lem7} For each $\alpha=1,\ldots,N$, $\Phi(x)\in C_0(\Omega)$
\begin{equation}\label{43}
F(\Phi B_{s\alpha})(\xi)\tilde\xi=0 \quad \forall\xi\in\R^n, \ \tilde\xi\in X_s.
\end{equation}
\end{lem}
\begin{pf}
Equality (\ref{43}) readily follows from the relation
$$
F(\Phi B_{s\alpha})(\xi)\tilde\xi=\int_{\R^n} e^{-2 \pi i\xi\cdot x}\Phi(x) B_{s\alpha}(x)\tilde\xi dx=0
$$
because $B_{s\alpha}(x)\tilde\xi=0$ for all $x\in\R^n$ by the definition of the subspace $X_s$.
\end{pf}

We introduce the set
\begin{eqnarray}\label{Lam}
\Lambda=\Lambda(x)=\Bigl\{ \lambda\in\C^N \ | \ \exists\eta\in\S: \nonumber\\
\sum_{\alpha=1}^N (iA_{s\alpha}\cdot\tilde\xi^s(\eta)-B_{s\alpha}\bar\xi^s(\eta)\cdot\bar\xi^s(\eta))\lambda_\alpha=0 \ \forall s=1,\ldots,m \ \Bigr\},
\end{eqnarray}
where $\tilde\xi^s(\eta)\in X_s$, $\bar\xi^s(\eta)\in X_s^\bot$ are such that $(\tilde\xi^s(\eta),\bar\xi^s(\eta))=p_{X_s}(\eta)$, and $p_{X_s}:\S\to S_{X_s}$ is the projection defined in section~\ref{Se2}. Let
$$
q(x,u)=Q(x)u\cdot u=\sum_{\alpha,\beta=1}^N q_{\alpha\beta}(x)u_\alpha \overline{u_\beta}
$$
be an Hermitian form with coefficients $q_{\alpha\beta}(x)\in C(\Omega)$.

Suppose that the sequence $q(x,u_r)\rightharpoonup v$ as $r\to\infty$ weakly in
$M_{loc}(\Omega)$. The following theorem is analogous to Theorems~\ref{Th1},\ref{Th3}.

\begin{thm}\label{Th4} If $q(x,\lambda)\ge 0$ for all $\lambda\in\Lambda(x)$,
$x\in\Omega$, then ${q(x,u(x))\le v}$.
\end{thm}
\begin{pf}
We fix $s\in\overline{1,m}$, and observe that in view of (\ref{42}) for each $\Phi(x)\in C_0^\infty(\Omega)$
the distributions
\begin{eqnarray}\label{44}
\sum_{\alpha=1}^N\sum_{k=1}^n\partial_{x_k}\left(a_{s\alpha k}\Phi u_{\alpha r}-2\sum_{l=1}^n b_{s\alpha kl}\Phi_{x_l} u_{\alpha r}\right)+
\sum_{\alpha=1}^N\sum_{k,l=1}^n \partial_{x_kx_l}(b_{s\alpha kl}\Phi u_{\alpha r})=\nonumber\\
\Phi\left(\sum_{\alpha=1}^N\sum_{k=1}^n \partial_{x_k}(a_{s\alpha k}u_{\alpha r})
+\sum_{\alpha=1}^N\sum_{k,l=1}^n \partial_{x_kx_l}(b_{s\alpha kl}u_{\alpha r})\right)+\nonumber\\ \sum_{\alpha=1}^N\sum_{k=1}^n a_{s\alpha k}\Phi_{x_k} u_{\alpha r}-\sum_{\alpha=1}^N\sum_{k,l=1}^n b_{s\alpha kl}\Phi_{x_kx_l} u_{\alpha r}
\end{eqnarray}
are pre-compact in the Sobolev space $H^{-1}(\R^n)\doteq W^{-1}_2(\R^n)$. Relation (\ref{44}) implies that the sequence
\begin{eqnarray}\label{45}
(1+|\xi|^2)^{-\frac{1}{2}}\left(\sum_{\alpha=1}^N\sum_{k=1}^n 2\pi i\xi_k \left(F(a_{s\alpha k}\Phi u_{\alpha r})(\xi)-2\sum_{l=1}^n F(b_{s\alpha kl}\Phi_{x_l} u_{\alpha r})(\xi)\right)\right. \nonumber\\ \left.
-\sum_{\alpha=1}^N\sum_{k,l=1}^n 4\pi^2\xi_k\xi_l F(b_{s\alpha kl}\Phi u_{\alpha r})(\xi)\right),
\end{eqnarray}
$r\in\N$, is compact in $L^2(\R^n)$. Denote $\tilde\xi=P_1\xi$, $\bar\xi=P_2\xi$, where $P_1,P_2$ are the orthogonal projections onto the subspaces $X_s$, $X_s^\bot$, respectively. Multiplying (\ref{45}) by the bounded function $(1+|\xi|^2)^{\frac{1}{2}}(1+|\tilde\xi|^2+|\bar\xi|^4)^{-\frac{1}{2}}$, we obtain that the sequence
\begin{eqnarray}\label{46}
(1+|\tilde\xi|^2+|\bar\xi|^4)^{-\frac{1}{2}}\sum_{\alpha=1}^N\bigl( 2\pi i (F(A_{s\alpha}\Phi u_{\alpha r})(\xi)-2F(B_{s\alpha}\nabla\Phi u_{\alpha r})(\xi))\cdot\xi\nonumber\\
- 4\pi^2F(B_{s\alpha}\Phi u_{\alpha r})(\xi)\xi\cdot\xi\bigr),
\end{eqnarray}
$r\in\N$, is compact in $L^2(\R^n)$.

By Lemma~\ref{Lem7} and symmetricity of the matrix $F(B_{s\alpha}\Phi u_{\alpha r})(\xi)$ we find that
\begin{eqnarray}\label{47}
F(B_{s\alpha}\Phi u_{\alpha r})(\xi)\xi\cdot\xi=F(B_{s\alpha}\Phi u_{\alpha r})(\xi)\bar\xi\cdot\bar\xi+
\nonumber\\ 2F(B_{s\alpha}\Phi u_{\alpha r})(\xi)\tilde\xi\cdot\bar\xi +F(B_{s\alpha}\Phi u_{\alpha r})(\xi)\tilde\xi\cdot\tilde\xi=F(B_{s\alpha}\Phi u_{\alpha r})(\xi)\bar\xi\cdot\bar\xi, \\
\label{48}
F(B_{s\alpha}\nabla\Phi u_{\alpha r})(\xi)\cdot\xi=F(B_{s\alpha}\nabla\Phi u_{\alpha r})(\xi)\cdot(\bar\xi+\tilde\xi)=F(B_{s\alpha}\nabla\Phi u_{\alpha r})(\xi)\cdot\bar\xi.
\end{eqnarray}
Notice also that the sequences
\begin{eqnarray}\label{49}
(1+|\tilde\xi|^2+|\bar\xi|^4)^{-\frac{1}{2}}F(A_{s\alpha}\Phi u_{\alpha r})(\xi)\cdot\bar\xi, \nonumber\\
(1+|\tilde\xi|^2+|\bar\xi|^4)^{-\frac{1}{2}}F(B_{s\alpha}\nabla\Phi u_{\alpha r})(\xi)\cdot\bar\xi, \ r\in\N, \mbox{ are compact in } L^2(\R^n),
\end{eqnarray}
since the functions $(1+|\tilde\xi|^2+|\bar\xi|^4)^{-\frac{1}{2}}\bar\xi_k$, $k=1,\ldots,n$ lay in the ideal $A_0$.
It now follows from (\ref{46})--(\ref{49}) that  the sequence of distributions
$$
l^s_r=(1+|\tilde\xi|^2+|\bar\xi|^4)^{-\frac{1}{2}}\sum_{\alpha=1}^N\bigl( 2\pi i F(A_{s\alpha}\Phi u_{\alpha r})(\xi)\cdot\tilde\xi
- 4\pi^2F(B_{s\alpha}\Phi u_{\alpha r})(\xi)\bar\xi\cdot\bar\xi\bigr),
$$
$r\in\N$, is compact in $L^2(\R^n)$.
The distributions $l^s_r$ can be represented as
\begin{equation}
l^s_r=\sum_{\alpha=1}^N \left( \sum_{k=1}^n 2\pi i p_{s\alpha k}(\xi)F(a_{s\alpha k}\Phi u_{\alpha r})(\xi)-\sum_{k,l=1}^n
4\pi^2q_{s\alpha kl}(\xi)F(b_{s\alpha kl}\Phi u_{\alpha r})(\xi)\right),
\end{equation}
where
\begin{eqnarray*}
p_{s\alpha k}(\xi)=(1+|\tilde\xi|^2+|\bar\xi|^4)^{-\frac{1}{2}}\tilde\xi_k\equiv (|\tilde\xi|^2+|\bar\xi|^4)^{-\frac{1}{2}}\tilde\xi_k \mod A_0, \\
q_{s\alpha kl}(\xi)=(1+|\tilde\xi|^2+|\bar\xi|^4)^{-\frac{1}{2}}\bar\xi_k\bar\xi_l\equiv (|\tilde\xi|^2+|\bar\xi|^4)^{-\frac{1}{2}}\bar\xi_k\bar\xi_l \mod A_0.
\end{eqnarray*}
In particular, we see that $p_{s\alpha k}(\xi),q_{s\alpha kl}(\xi)\in A_{X_s}$.

Taking into account compactness of commutators $[A_\psi,B_\phi]$ in $L^2(\R^n)$, where $(\psi,\phi)=(p_{s\alpha k}(\xi),\chi(x) a_{s\alpha k}(x))$, $(\psi,\phi)=(q_{s\alpha kl}(\xi),\chi(x)b_{s\alpha kl}(x))$, and $\chi(x)\in C_0(\Omega)$ is a function such that $\chi(x)\Phi(x)=\Phi(x)$, we find that
the sequence
$$
\chi(x)\sum_{\alpha=1}^N \left( 2\pi i\sum_{k=1}^n a_{s\alpha k}(x)p_{s\alpha k}(\partial/\partial x)-4\pi^2 \sum_{k,l=1}^n b_{s\alpha kl}(x) q_{s\alpha kl}(\partial/\partial x)\right)(\Phi u_{\alpha r}),
$$
$r\in\N$, is compact in $L^2(\R^n)$, that is, the sequence
$$
\sum_{\alpha=1}^N \left( 2\pi i\sum_{k=1}^n a_{s\alpha k}(x)p_{s\alpha k}(\partial/\partial x)-4\pi^2 \sum_{k,l=1}^n b_{s\alpha kl}(x) q_{s\alpha kl}(\partial/\partial x)\right)(\Phi u_{\alpha r}),
$$
$r\in\N$, is compact in $L^2_{loc}(\Omega)$. Here $p_{s\alpha k}(\partial/\partial x)$, $q_{s\alpha kl}(\partial/\partial x)$ are pseudodifferential operators with symbols
$$
(|\tilde\xi|^2+|\bar\xi|^4)^{-\frac{1}{2}}\tilde\xi_k, \ (|\tilde\xi|^2+|\bar\xi|^4)^{-\frac{1}{2}}\bar\xi_k\bar\xi_l
$$
laying in $A_{X_s}\subset\A$. Since $s=1,\ldots,m$, $\Phi(x)\in C_0^\infty(\Omega)$ are arbitrary, we see that our sequence $u_r$ satisfies constraints of the kind (\ref{38}).
Since
$$
\widehat{p_{s\alpha k}}(\eta)=\tilde\xi^s_k(\eta), \ \widehat{q_{s\alpha kl}}(\eta)=\bar\xi^s_k(\eta)\bar\xi^s_l(\eta),
$$
the set $\Lambda=\Lambda(x)$ corresponding to these constraints is
\begin{eqnarray*}
\Lambda=\Lambda(x)=\Bigl\{ \lambda\in\C^N \ | \ \exists\eta\in\S: \\
\sum_{\alpha=1}^N (2\pi iA_{s\alpha}\cdot\tilde\xi^s(\eta)-4\pi^2B_{s\alpha}\bar\xi^s(\eta)\cdot\bar\xi^s(\eta))\lambda_\alpha=0 \ \forall s=1,\ldots,m \ \Bigr\}.
\end{eqnarray*}
Since, in accordance with Remark~\ref{Rem1},  $\tilde\xi^s(t\eta)=a(t,\eta)\tilde\xi^s(t\eta)$, $\bar\xi^s(t\eta)=b(t,\eta)\bar\xi^s(t\eta)$, $b^2(t,\eta)=ta(t,\eta)$, then after the transformation $\eta=(2\pi)^{-1}\eta$ the set $\Lambda$ will coincide with (\ref{Lam}). Then the assertion of Theorem~\ref{Th4} readily follows from Theorem~\ref{Th3}. The proof is complete.
\end{pf}

\subsection{One example}

Let us consider the sequence $u_r=(u_{r1},u_{r2},u_{r3})\in L^2_{loc}(\Omega,\C^3)$, $\Omega\subset\R^3$ weakly convergent to $u=(u_1,u_2,u_3)$ such that the sequences
$$
i\left(\partial_{x_3} u_{r2}-\partial_{x_2} u_{r3}\right)+\partial_{x_1}^2 u_{r1}; \quad
i\left(\partial_{x_1} u_{r3}-\partial_{x_3} u_{r1}\right)+\partial_{x_2}^2 u_{r2}; \quad
i\left(\partial_{x_2} u_{r1}-\partial_{x_1} u_{r2}\right)+\partial_{x_3}^2 u_{r3}
$$
are pre-compact in $H^{-1}_{loc}(\Omega)$.

\begin{thm}\label{Th5} For every pair $(k,l)$, $1\le k<l\le 3$ we have
$$ u_{rk}\overline{u_{rl}}\mathop{\rightharpoonup}_{r\to\infty} u_k\overline{u_l}.$$
\end{thm}
\begin{pf}
In the notations of Theorem~\ref{Th4} we find that $X_i=\{\xi\in\R^3: \xi_i=0\}$, $i=1,2,3$,
while the set $\Lambda$ is determined by the relations
\begin{eqnarray}\label{50}
\lambda_2\tilde\xi^1_3(\eta)-\lambda_3\tilde\xi^1_2(\eta)+\lambda_1(\bar\xi^1_1(\eta))^2=
\lambda_3\tilde\xi^2_1(\eta)-\lambda_1\tilde\xi^2_3(\eta)+\lambda_2(\bar\xi^2_2(\eta))^2\nonumber\\
=\lambda_1\tilde\xi^3_2(\eta)-\lambda_2\tilde\xi^3_1(\eta)+\lambda_3(\bar\xi^3_3(\eta))^2=0
\end{eqnarray}
for some $\eta\in\S$. For $\gamma\in\C$ we introduce the Hermitian form
$$Q_\gamma(\lambda)=\mathop{Re} \gamma u_k\overline {u_l}=\frac{\gamma}{2}u_k\overline {u_l}+
\frac{\bar\gamma}{2}u_l\overline {u_k}.
$$
Let $\lambda\in\Lambda$. Then there exists $\eta\in\S$ such that (\ref{50}) holds.
Observe that the space $\tilde X$ from Proposition~\ref{Pro3} may be included at most in two subspaces $X_i$.
If the set $I=\{ \alpha\in\overline{1,3} : \tilde X\not\subset X_i \}$ contains two different indexes $j,k$, then
$\tilde\xi^j(\eta)=\tilde\xi^k(\eta)=0$, $|\bar\xi^j(\eta)|=|\bar\xi^k(\eta)|=1$ by Proposition~\ref{Pro3} and it follows from (\ref{50}) that $\lambda_k=\lambda_j=0 \Rightarrow Q_\gamma(\lambda)=0$.

In the remaining case there exists only one index $j$ such that $\tilde X\not\subset X_j$.
For definiteness, we assume that $j=1$. Then again $\tilde\xi^1(\eta)=0$, $\bar\xi^1(\eta)\not=0$, which imply that $\lambda_1=0$. It is clear that $\tilde X=X_2\cap X_3$. By Proposition~\ref{Pro3} we find that $\tilde\xi^2(\eta)=(a,0,0)$, $\tilde\xi^3(\eta)=(b,0,0)$, and $ab>0$.
Let $\bar\xi^2(\eta)=(0,p,0)$, $\bar\xi^3(\eta)=(0,0,q)$.
By (\ref{50})  $$a\lambda_3+p^2\lambda_2=-b\lambda_2+q^2\lambda_3=0.$$ Since the determinant of this system $\Delta=p^2q^2+ab>0$ we conclude $\lambda_2=\lambda_3=0$. Thus, $\lambda=0$ and $Q_\gamma(\lambda)=0$. By Theorem~\ref{Th4} we see that $Q_\gamma(u_r)\rightharpoonup Q_\gamma(u)$.
Therefore,
$$
u_{rk}\overline{u_{rl}}=Q_1(u_r)-iQ_i(u_r)\mathop{\rightharpoonup}_{r\to\infty}Q_1(u)-iQ_i(u)=u_k\overline{u_l},
$$
as was to be proved.
\end{pf}
Observe that in the notations of Theorem~\ref{Th1} the set $\Lambda=\{\lambda\in\R^3 | \lambda_1\lambda_2\lambda_3=0\}$ and this theorem does not allow to derive the statement of Theorem~\ref{Th5}.


\end{document}